\documentclass{amsart}
\usepackage{color,graphicx,wrapfig}
\newcommand{\non}{\nonumber}
\title[]{On the singularities of a free boundary\\
through Fourier expansion}
\author[J. Andersson]{John Andersson}
\address{Mathematics Institute,
University of Warwick
Coventry CV4 7AL, UK}
\email{j.e.andersson@warwick.ac.uk}
\author[H. Shahgholian ]{Henrik Shahgholian}
\address{Department of Mathematics, Royal Institute of Technology,
100~44  Stockholm, Sweden}
\email{henriksh@math.kth.se}
\urladdr{http://www.math.kth.se/~henriksh/}
\author[G.S. Weiss]{Georg S. Weiss}
\address{Graduate School of Mathematical Sciences,
University of Tokyo, 3-8-1 Komaba, Meguro-ku, Tokyo-to, 153-8914 Japan,}
\email{gw@ms.u-tokyo.ac.jp}
\urladdr{http://www.ms.u-tokyo.ac.jp/~gw/}
\thanks{$2000$ {\it Mathematics Subject Classification.\/} Primary
35R35, Secondary 35B40, 35J60.}
\thanks{H. Shahgholian has been supported in part by
the Swedish Research Council.
G.S. Weiss has been partially supported by the Grant-in-Aid
21540211 of the Japanese Ministry of Education, Culture, Sports, Science and Technology.
He also thanks the Knut och Alice Wallenberg foundation for a visiting appointment to
KTH.
Both J. Andersson and G.S. Weiss thank the G\"oran Gustafsson Foundation
for visiting appointments to
KTH}
\date{}
\def\dh{{\> d\mathcal{H}^{n-1}}}
\def\dh1{{\> d\mathcal{H}^1}}
\def\dh2{{\> d\mathcal{H}^2}}
\def\P{{\mathbb{P}_2}}
\def\R{{\mathbb{R}}}
\def\L{{\mathcal{L}}}

\def\dist{\hbox{\rm dist}}

\def\x{{\mathbf x}}
\def\y{{\mathbf y}}
\theoremstyle{plain}
\newtheorem{theorem}{Theorem}[section]

\newtheorem{lem}[theorem]{Lemma}
\newtheorem{prop}[theorem]{Proposition}
\newtheorem{cor}[theorem]{Corollary}
\newtheorem{thm}[theorem]{Theorem}
\theoremstyle{definition}
\newtheorem{definition}[theorem]{Definition}
\theoremstyle{definition}
\newtheorem{remark}[theorem]{Remark}
\theoremstyle{example}

\theoremstyle{example}
\newtheorem*{discussion}{Discussion}
\numberwithin{equation}{section}
\begin{document}
\maketitle
\begin{abstract}
In this paper we are concerned with singular points of solutions to the {\it unstable} free boundary problem
$$
\Delta u = - \chi_{\{u>0\}} \qquad \hbox{ in } B_1.
$$
The problem arises in  applications such as solid combustion, composite membranes, climatology and fluid dynamics.

It is known that solutions to the above problem may exhibit
singularities ---that is  points at which the second derivatives of the solution are unbounded---
as well as degenerate points. This causes breakdown of by-now classical techniques.
Here we introduce new ideas based on Fourier expansion of the nonlinearity $\chi_{\{u>0\}} $.

The method turns out to have enough momentum to accomplish a complete description of the structure of the singular set
in ${\mathbb R}^3$.

A surprising fact in ${\mathbb R}^3$ is that although
$$
\frac{u(r\x)}{\sup_{B_1}|u(r\x)|}
$$
can converge at singularities to each of the harmonic polynomials
$$ xy, {x^2+y^2\over 2}-z^2 \textrm{ and } z^2-{x^2+y^2\over 2},$$
it may {\em not} converge to any of the non-axially-symmetric harmonic polynomials
$\alpha\left((1+ \delta )x^2 +(1- \delta )y^2 - 2z^2\right)$
with $\delta\ne 1/2$.

We also prove the existence of stable
singularities in ${\mathbb R}^3$.
\end{abstract}
\tableofcontents
\section{Introduction}

We investigate the singular points of solutions of
the unstable free boundary problem 
\begin{equation}\label{main}
\Delta u=-\chi_{\{u>0\}} \textrm{ in }B_1,
\end{equation}
arising in solid combustion (see the references in \cite{monneauweiss}),
the composite membrane problem (\cite{chanillo1},
\cite{chanillo2}, \cite{blank}, \cite{shahgholian}, \cite{chanillokenig}, \cite{chanillokenigto}),
climatology (\cite{diaz})
and fluid dynamics (\cite{struwe}).

The minus sign on the right-hand side drastically changes the
problem from the well-known obstacle problem 
(see for example \cite{ca1}, \cite{ca2} and \cite{ca3})
into an
unstable problem exhibiting non-uniqueness, bifurcations,
unbounded second derivatives and more. Let us
describe some of the known results. 

From standard elliptic regularity theory it follows that
if $u$ is a solution to (\ref{main}) then $u\in C^{1,\alpha}$ for all $\alpha<1$.
However, in contrast to the well-known obstacle problem $\Delta u=\chi_{\{u>0\}}$,
the solutions to (\ref{main}) are not $C^{1,1}$ in general. The existence of
non-regular solutions was first shown in \cite{anderssonweiss}.

For convenience let us denote the set of singular points by
$$S^u=\{\x :  u\notin C^{1,1}(B_r(\x)) \textrm{ for any } r>0\}.$$
As observed in \cite{monneauweiss}, $\{ u=0\}$ is analytic 
and $u\in C^{1,1}$
in 
a neighborhood of each $\x\in \{ u=0\}\cap \{ \nabla u\ne 0\}$.
Thus the set of singular points
is contained in the set where both $u$ and $|\nabla u|$ vanish.

We may expect that for $\x^0\in S^u$ the blow-up
\begin{equation}\label{genblowup}
\lim_{r\to 0}\frac{u(r\x+\x^0)}{\sup_{B_r(\x^0)}|u|},
\end{equation}
should give us some information about the singular set.
It was shown in \cite{monneauweiss} (see also Proposition \ref{fixedcenter} below)
that at a singular point $\x^0$
$$
\lim_{j\to \infty}\frac{u(r_j\x+\x^0)}{\sup_{B_{r_j}(\x^0)}|u|}=p(\x),
$$
where $p$ is a second order homogeneous harmonic polynomial.
This raises several questions.
\begin{enumerate}
\item Does $p$ depend on the choice of the sequence $r_j\to 0$? 
\item Does every second order homogeneous harmonic polynomial $p$
appear as limit?
\item Is there any partial regularity of the singular set?
\item Do energy minimising singularities exist?
\end{enumerate}
Concerning uniqueness of blow-up limits it has been shown in
\cite{monneauweiss} that
in two
dimensions the free boundary of the {\em minimal solution}
close to points where the
second derivative is unbounded, consists of four Lip\-schitz graphs meeting
at right angles.
In \cite{cmp} this fact has been extended to any solution
in two dimensions, proving also uniformity and quantitative
estimates by methods closely related to those in the present
paper.
Concerning question (iv) it has been proved in \cite{monneauweiss} 
that the singularity in two dimensions is {\em unstable}
in the sense that the second variation of the energy is
negative.
As to stability of higher dimensional singularities
there is a gap in the proof of \cite{monneauweiss}
which has been pointed out by Carlos Kenig-Sagun Chanillo-Tung To (\cite{chanillokenigto}).
In the present paper we will prove the following main results
which among other things close the gap in \cite{monneauweiss}
by showing that the cut-off dimension concerning this
problem is $3$, that is, there exist stable singularities
in three dimensions.

{\bf Main results:}
{\em\begin{enumerate}
\item Existence of a true three-dimensional singularity (Corollary \ref{exist}).
\item Axial symmetry of blow-up limits in three dimensions (Theorem \ref{symm}).
\item Unique tangent cones at true three-dimensional singularities in ${\mathbb R}^3$ (Theorem \ref{Z_1}).
\item Unique tangent cones at unstable codimension two singularities in ${\mathbb R}^3$ (Theorem \ref{uniqqesinglines}).
\item Stability of true three-dimensional singularities in ${\mathbb R}^3$ (Theorem \ref{stable}).
\item Regularity of the singular set in three dimensions (Section \ref{regularityofthesingularset}). 
\end{enumerate}
}
\begin{discussion}
In contrast to the analysis of singularities for minimisers or stable
solutions, where there are many methods available, there are few
results on unique tangent cones at {\em unstable} singularities.
Even the Lojasiewicz inequality approach (see for example \cite{simon})
would be hard to realize in our problem due to the
lack of a suitable local Lyapunov functional; we do have a monotonicity
formula playing the role of a local Lyapunov functional, but
as it turns out it has the wrong scaling to be used at the
unstable singularities of ``supercharacteristic growth''. 

The natural approach would be to study blow-up limits in order to
analyze the singularities.
Unfortunately
the blow-up sequence in
(\ref{genblowup}) does not provide enough information of the solution as
the nonlinearity of equation (\ref{main}) vanishes in the limit. To
preserve some information of the nonlinearity we will instead, 
in Section \ref{genback},
consider
\begin{equation}\label{typetwoblowup}
\frac{u(r_j\x+\x^0)}{r_j^2}-\Pi(u,r_j,\x^0),
\end{equation}
where $\Pi(u,r_j,\x^0)$ is the projection of $u(r_j\x+\x^0)/r_j^2$ in $B_1$ onto
the homogeneous harmonic second order polynomials (see Definition
\ref{projection}).

It can be shown that if
\begin{align*}
&\lim_{j\to \infty}\frac{u(r_j\x+\x^0)}{\sup_{B_{r_j}(\x^0)}|u|}=p(\x),
\intertext{then}
&\lim_{j\to \infty}\left(\frac{u(r_j\x+\x^0)}{r_j^2}-\Pi(u,r_j,\x^0)\right)=Z_p,
\end{align*}
where $Z_p$ is a solution of
$$
\Delta Z_p=-\chi_{\{p>0\}}.
$$
Next we notice that, at each singular point $\x^0$,
$$
\lim_{j\to \infty}\frac{\Pi(u,r_j,\x^0)}{\sup_{B_1}|\Pi(u,r_j,\x^0)|}=
\lim_{j\to \infty}\frac{u(r_j\x+\x^0)}{\sup_{B_{r_j}(\x^0)}|u|}.
$$
So in order to prove uniqueness of $p$ it is sufficient to
control how $\Pi(u,r,\x^0)$ changes when $r$ varies. More precisely
we would want to estimate
\begin{equation}\label{pidiff}
\Big|\frac{\Pi(u,r,\x^0)}{\sup_{B_1}|\Pi(u,r,\x^0)|}-
\frac{\Pi(u,r/2,\x^0)}{\sup_{B_1}|\Pi(u,r/2,\x^0)|}\Big|.
\end{equation}

Our method of proof is based on the observation that
$u(r\x+\x^0)\approx \tau_r p_r +Z_{p_r}$ in $B_r$, where
$p_r$ is a second order harmonic polynomial of norm $1$. 
It follows that $\Pi(u,r/2,\x^0)\approx \Pi(\tau_r p_r+Z_{p_r},1/2,0)=
\Pi(\tau_rp_r,1/2,0)+\Pi(Z_{p_r},1/2,0)=\tau_rp_r +\Pi(Z_{p_r},1/2,0)$
(cf. Section
\ref{prelanal}).
Therefore it is essential to control $\Pi(Z_{p_r},\cdot)$
in order to estimate (\ref{pidiff}).
This control will be achieved by means of an explicit calculation
of the Fourier coefficients of $Z_{p_r}$.
\end{discussion}

\noindent{\bf Plan of the paper.}\\ 
In Section \ref{genback}
we will remind ourselves of results and definitions
of \cite{cpde}, \cite{monneauweiss},
\cite{anderssonweiss} and \cite{cmp}
that are relevant to the present paper.

In Section \ref{Fouriersec} we use techniques
developed in \cite{KarpMargulis} 
based on Fourier coefficients
to analyze $Z_{p}$. We also explicitly
calculate $Z_{p}$ when $p=2xz$ or $p=\pm\big((x^2+y^2)/2-z^2 \big)$.
Using these calculations we are able to show in Section \ref{growth} that
in three dimensions and for
small $r$
$$\sup_{B_r}|u(\x+\x^0)|\ge c r^2|\log r|$$
(Corollary \ref{loggrowth}).

Based on this estimate on the growth of $u$ we prove
in Section \ref{truethree} existence of a true three-dimensional
singularity. 

Section \ref{prelanal}
provides estimates on 
$u-\Pi(u)-Z_{\Pi(u)}$
which we combine in Section \ref{SecClass}
carefully with 
the above analysis of
$Z_p$ to 
show axial symmetry of blow-up limits in three dimensions.
This is a remarkable symmetrization effect in view of the fact
that there are of course second order
homogeneous harmonic polynomials that are {\em not} axially symmetric.

In Sections \ref{uniquesingpoints} and \ref{uniquesinglines} we prove
---once more carefully using the information gained on the
Fourier coefficients---
uniqueness of the blow-up limits at
singular points in ${\mathbb R}^3$. 
Based on the asymptotics in Section \ref{uniquesingpoints}
we are able to show in Section \ref{cone}
that true three-dimensional singularities are stable.
In Section \ref{regularityofthesingularset}
we use standard techniques to show that in the three-dimensional case
the singular set may be decomposed into a countable set of
isolated points and a component that
is locally
contained in a $C^1$-curve.
In the Appendix
we have gathered technical calculations 
which may well be considered to be the core of our paper.

\section{Notation}
Throughout this article ${\mathbb R}^n$ will be equipped with the Euclidean
inner product $\x\cdot \y$ and the induced norm $\vert \x \vert\> .$
Moreover $A:B=\sum_{i,j=1}^n a_{ij} b_{ij}$ shall denote the inner product
of two $(n,n)$ matrices.
We will use the set $\mathcal{Q}$
of all orthogonal matrices in $ {\mathbb R}^n$.
We define
$B_r(\x^0)$ as the open $n$-dimensional ball of center
$\x^0\> ,$ radius $r$ and volume $r^n\> \omega_n$,
and $B^+_r(\x^0):= \{ \x\in B_r(\x^0) : x_n>0\},
B^-_r(\x^0):= \{ \x\in B_r(\x^0) : x_n<0\}$.
When not specified, $\x^0$ is assumed to be $0$.
We shall often use
abbreviations for inverse images like $\{u>0\} :=
\{\x\in \Omega\> : \> u(\x)>0\}\> , \> \{x_n>0\} :=
\{\x \in {\mathbb R}^n \> : \> x_n > 0\}$ etc.
and occasionally
we shall employ the decomposition $\x=(x_1,\dots,x_n)$ of a vector $\x\in {\mathbb R}^n\> .$
We will use the $n$-dimensional Lebesgue measure
${\mathcal L}^n$ and the $k$-dimensional Hausdorff measure
${\mathcal H}^k$.
When considering a set $A\> ,$ $\chi_A$ shall stand for
the characteristic function of $A\> ,$
while
$\nu$ shall typically denote the outward
normal to a given boundary.
We will use Landau's symbols as signed variables.
For example $-o(1)$ will mean a negative quantity that
turns to zero.
By $\P$ we will denote the space of second order homogeneous
harmonic polynomials in $\R^n$.
We shall also use the projection
$\Pi$ onto $\P$ as well as the norm $\tau$ of $\Pi(v)$,
both defined in Definition \ref{projection},
as well as the parametrization parameters
$\delta^A(v)$ and $\delta^B(v)$ defined in Definition \ref{delta}.
Last, we shall use for $p\in \P$ the Newtonian potential $Z_p$, i.e.
the unique solution of 
\begin{align}
&\Delta Z_p = -\chi_{\{ p>0\}} \textrm{ in } \R^n,\non\\
&Z_p(0)=|\nabla Z_p(0)|=0,\non\\
&\lim_{|\x|\to \infty} \frac{Z(\x)}{|\x|^3}=0\textrm{ and}\non\\
&\Pi(Z_p,1)=0\textrm{ (cf. \cite{cmp} Section 4).}\non
\end{align}
\section{General Background}\label{genback}

In this section we will gather some results from \cite{monneauweiss} and
\cite{anderssonweiss}, and describe some compactness properties
of blow-ups of solutions. First we will remind ourselves of the monotonicity
formula proved in \cite{cpde}.
The roots of those monotonicity formulas
are harmonic mappings (\cite{schoen}, \cite{price}) and blow-up (\cite{pacard}).

\begin{thm}\label{mon}
Suppose that $u$ is a solution of (\ref{main}) in $\Omega$
and that $B_\delta(\x^0)\subset \Omega\> .$
Then for all $0<\rho<\sigma<\delta$
the function
\[ \Phi^u_{\x^0}(r) := r^{-n-2} \int_{B_r(\x^0)} \left(
{\vert \nabla u \vert}^2 \> -\> 2\max(u,0)
\right)\]\[
- \; 2 \> r^{-n-3}\>  \int_{\partial B_r(\x^0)}
u^2 \> d{\mathcal H}^{n-1}\; ,\]
defined in $(0,\delta)\> ,$ satisfies the monotonicity formula
\[ \Phi^u_{\x^0}(\sigma)\> -\> \Phi^u_{\x^0}(\rho) \; = \;
\int_\rho^\sigma r^{-n-2}\;
\int_{\partial B_r(\x^0)} 2 \left(\nabla u \cdot \nu - 2 \>
{u \over r}\right)^2 \; d{\mathcal H}^{n-1} \> dr \; \ge 0 \; \; .\]
\end{thm}
This energy monotonicity is important since it helps us to distinguish
different points of the set $\{u(0)=|\nabla u(0)|=0\}$. In particular
we may according to the following Proposition define the singular set $S^u$ as
$$
S^u=\{\x\in B_1:\; u(\x)=|\nabla u(\x)|=0 \textrm{ and } \lim_{r\to
0}\Phi^u_{\x}(r)=-\infty \}.
$$
\begin{prop}[Proposition 5.1 in
\cite{monneauweiss}]\label{fixedcenter}
Let $u$ be a solution of (\ref{main}) in $\Omega$
and let us consider a point
$\x^0\in \Omega\cap \{ u=0\}\cap\{ \nabla u =0\}.$\\
(i) In the case $\Phi^u_{\x^0}(0+)=-\infty$,
$\lim_{r\to 0} r^{-3-n}\int_{\partial B_r(\x^0)} u^2 \> d{\mathcal H}^{n-1}
= +\infty$, and for
$T(\x^0,r) := \left(r^{1-n}\int_{\partial B_{r}(\x^0)} u^2\> d{\mathcal H}^{n-1}
\right)^{1\over 2},$
each limit of
\[ \frac{u(\x^0+r \x)}{T(\x^0,r)}\]
as $r\to 0$ belongs to $\P$.
\\
(ii) In the case $\Phi^u_{\x^0}(0+)\in (-\infty,0)$,
\[ u_r(\x) := \frac{u(\x^0+r \x)}{r^2}\]
is bounded in $W^{1,2}(B_1(0))$,
and each limit as $r\to 0$ is a homogeneous solution of degree $2$.\\
(iii) If $\Phi^u_{\x^0}(0+)=0$, then
\[ \frac{u(\x^0+r \x)}{r^2}\to 0\hbox{ in } W^{1,2}(B_1(0)) \hbox{ as } r\to 0\; .\]
\end{prop}
\begin{remark}\label{int}
In \cite[Lemma 5.2]{monneauweiss} it says that case (ii) in
Proposition \ref{fixedcenter} does not occur in $\mathbb{R}^2$. Unfortunately
the authors omitted the following homogeneous solution of second order
in $\mathbb{R}^2$:

Using polar coordinates $(r,\phi)$, let
$$
u(r,\phi)=\left\{
\begin{array}{ll}
-\frac{r^2}{2}\cos^2(\phi)+\frac{r^2}{4\sqrt{3}}\sin(2\phi) & \textrm{ when }
\phi\in (\pi/3,\pi/2), \\
-\frac{r^2}{8}\cos(2\phi)-\frac{r^2}{8\sqrt{3}}\sin(2\phi) & \textrm{ when }
\phi\in (-\pi/6,\pi/3), \\
-\frac{r^2}{2}\cos^2(\phi+2\pi/3)+\frac{r^2}{4\sqrt{3}}\sin(2\phi+2\pi/3) & \textrm{
when } \phi\in (-\pi/3,-\pi/6), \\
-\frac{r^2}{8}\cos(2\phi+2\pi/3)-\frac{r^2}{8\sqrt{3}}\sin(2\phi+2\pi/3) & \textrm{
when }\phi\in (-5\pi/6,-\pi/3), \\
-\frac{r^2}{2}\cos^2(\phi+4\pi/3)+\frac{r^2}{4\sqrt{3}}\sin(2\phi+4\pi/3) & \textrm{
when }\phi\in (-\pi,-5\pi/6), \\
-\frac{r^2}{8}\cos(2\phi+4\pi/3)-\frac{r^2}{8\sqrt{3}}\sin(2\phi+4\pi/3) & \textrm{
when } \phi\in (-\pi,-3\pi/2).
\end{array}\right.
$$
Then $u$ is a second order homogeneous solution to
equation (\ref{main}).

Let us show that up to rotations, $u$ is the unique non-trivial second
order
homogeneous solution to equation (\ref{main}) in $\mathbb{R}^2$.
Each cone in which $u$ is 
negative has to have an opening of exactly
$\pi/2$. Thus $u$ can be negative in at most $4$ different connected components.
However if $u$ is negative in four components then $u\le 0$ and
thus $\Delta u=0$. Since $u(0)=0$, $u\equiv 0$ by the strong maximum
principle. If $u$ is negative in only one component then $\Delta u=-1$
in a cone with opening $3\pi/2$ with zero boundary values on that cone, and
it is easy to see that such a $u$ is not homogeneous of second order.
If $u\le 0$ in two components then, after a rotation, $\Delta u=-\chi_{\{xy>0\}}$.
In ``$2$-dimensional Solutions'' (vi) in Section \ref{Fouriersec}, we will see that such a solution is
not homogeneous either. The only remaining possibility is that $u\le 0$ in three components, with
three components of $u> 0$ in between. As the gradient of $u$ is
continuous across the zero level set and $u$ is symmetric in each cone
where $u$ has a sign it follows that the opening of the cones where
$u> 0$ must equal each other. It follows that $u> 0$ in three
cones of opening $\pi/6$ where $\Delta u=-1$. Thus
$u$ is unique up to a rotation.
\end{remark}
In \cite{anderssonweiss} the authors have obtained
existence of solutions
in two dimensions exhibiting {\em cross-like singularities}
at which the second derivatives of the solution are
unbounded (case (i) of Proposition \ref{fixedcenter}), as well as degenerate
singularities at which
the solution decays to zero faster than any quadratic
polynomial (case (iii) of Proposition \ref{fixedcenter}):
\begin{thm}[Cross-shaped singularity, Corollary 4.2 in
\cite{anderssonweiss}]\label{cross}
There exists a solution $u$ of
\begin{displaymath}
\Delta u= -\chi_{\{u>0\}} \quad \textrm{in } B_1\subset \mathbb{R}^2
\end{displaymath}
that is {\bf not} of class $C^{1,1}$.
Each limit of
\[ \frac{u(r\x)}{T(0,r)}\]
as $r\to 0$ coincides after rotation with the function
$(x_1^2-x_2^2)/\Vert x_1^2-x_2^2\Vert_{L^2(\partial B_1(0))}$.
\end{thm}
\textsl{Remark on the Proof:} In \cite{anderssonweiss} the authors
show that one can construct a solution $u$ to (\ref{main}) with
$\lim_{r\to 0}\Phi^u_{0}(r)\le -M$ for any $M\ge 0$.
Then they use \cite[Lemma 5.2]{monneauweiss} that states that if
$\lim_{r\to 0}\Phi^u_{0}(r) < 0$ in $\mathbb{R}^2$ then we are in
case $(i)$ of Proposition \ref{fixedcenter}. As we pointed out in Remark
\ref{int}, Lemma 5.2 in \cite{monneauweiss} is not true.
The proof in \cite{anderssonweiss} however is easily fixed:
we only have to notice that
all second order homogeneous solutions of (\ref{main}) have fixed
energy $\Phi^u_0=m_0$ which follows from the uniqueness in Remark \ref{int}.
Thus if we choose the constant $M$ large enough we can exclude the possibility
that we are in case $(ii)$ or $(iii)$ of Proposition \ref{fixedcenter}
and the theorem follows.\qed

The proof of the previous theorem can be adjusted to
construct other kinds of singular points
(see Corollary \ref{exist}).

\begin{definition}\label{projection}
By $\Pi(u,r,\x^0)$ we will denote the projection operator onto $\P$ defined as follows:
$\Pi(u,r,\x^0)=\tau p$, where $\tau\in \mathbb{R}^+$ and $p\in \P$ satisfies
$\sup_{B_1}|p|=1$ as well as
$$
\inf_{h\in \P}\int_{B_1(0)}\Big|\frac{D^2u(r\x+\x^0)}{r^2}-D^2 h \Big|^2=
\int_{B_1(0)}\Big|\frac{D^2u(r\x+\x^0)}{r^2}-\tau D^2 p \Big|^2.
$$

We will often write $\Pi(u,r)$ when $\x^0$ is either the origin or given by
the context. At times we will also denote $\tau_r=\sup_{B_1}|\Pi(u,r)|$
and $p_r=\Pi(u,r)/\tau_r$.
\end{definition}
The following Lemma justifies the previous Definition.

\begin{lem}\label{proj}
The following four statements hold.\\
(i) For each $v\in W^{2,2}(B_1)$ the minimizer of
Definition \ref{projection} exists and is unique.
Thus $\Pi:W^{2,2}(B_1)\times (0,s)\times B_{1-s}\to P$ is well-defined.
\\
(ii) $\Pi$ is a linear operator.\\
(iii) If $h\in W^{2,2}(B_1)$ is harmonic in $B_1$
then $\Pi(h(\x),s)=\Pi(h(\x),r)$ for all $s,r\in (0,1)$.\\
(iv) For every $v,w\in W^{2,2}(B_1)$,
$$\sup_{B_1}|\Pi(v+w,r)|\le \sup_{B_1}|\Pi(v,r)|\> +\> \sup_{B_1}|\Pi(w,r)|.$$
\end{lem}
\proof
The first and second statement follow from the projection theorem
with respect to the $L^2(B_1;\mathbb{R}^{n^2})$-inner product and the linear subspace
$\{ f\in  L^2(B_1;\mathbb{R}^{n^2})\> : \> f(\x)
\textrm{ is symmetric, constant and trace}(f(\x)) = 0\}$.
\\
Writing $h$ as the sum of homogeneous harmonic polynomials $h_j$
that are orthogonal to each other with respect to
$$(v,w) := \int_{B_1} \sum_{i,j=1}^n \partial_{ij} v \partial_{ij} w,$$
we see that $\Pi(h_j)=0$ for all $j$ such that the degree of $h_j$
is different from $2$, implying the third statement.
\\
The last statement follows from the linearity of $\Pi$ and
the triangle inequality in $L^2(B_1;\mathbb{R}^{n^2})$.
\qed

Next we mention that solutions to (\ref{main}) have second derivatives
in $BMO$. This has been proved in \cite[Lemma 5.1]{cmp}
using standard facts of harmonic analysis.
% In Appendix B we supply a new and different proof based
% on a blow-up argument. 
\begin{prop}[cf. \hbox{\cite[Lemma 5.1]{cmp}}]\label{boundfordiff}
Let $u$ be a solution to (\ref{main}) in $B_1$ such that $\sup_{B_1}|u|\le M$ and
$u(0)=|\nabla u(0)|=0$. Then 
$$
\sup_{B_1}\Big| \frac{u(r\x)}{r^2}-\Pi(u,r) \Big|\le C_0
\textrm{ for every } r\le {1\over 2}
$$
where the constant $C_0$ depends only on $M$ and $n$.

Furthermore, for each $\alpha<1$
$$
\big{\|} \frac{u(r\x)}{r^2}-\Pi(u,r) \big{\|}_{C^{1,\alpha}(\overline{B_1})}\le C_\alpha(M,n),
$$
and for each $p< \infty$
$$
\big{\|} \frac{u(r\x)}{r^2}-\Pi(u,r) \big{\|}_{W^{2,p}(B_1)}\le C_p(M,n).
$$
\end{prop}
\section{Fourier Series Expansions of Global Solutions}\label{Fouriersec}

In this section we will remind ourselves of the work of L. Karp and A.S.
Margulis \cite{KarpMargulis}. In particular,
Theorem 3.1 and Proposition 3.2 in \cite{KarpMargulis}, summarized in
the next theorem, will be of importance to us.

\begin{thm}[\cite{KarpMargulis}]\label{KarpM}
Let $\sigma\in L^\infty(\mathbb{R}^n)$ be homogeneous of
zeroth order, that is $\sigma(\x)=\sigma(r\x)$ for all $r>0$. Assume that
$\sigma$ has the Fourier series expansion
$$
\sigma(\x)=\sum_{i=0}^\infty a_i \sigma_i,
$$
on the unit sphere, where $\sigma_i$ is a homogeneous harmonic polynomial of order
$i$.

Moreover assume that
$\Delta u= \sigma$ and that
$u(0)=|\nabla u(0)|=\lim_{\x\to \infty}u(\x)/|\x|^3=0$. Then
\begin{align*}
&u(\x)=q(\x)\log|\x|+|\x|^2\phi(\x),
\intertext{where}
&q=\frac{a_2}{n+2}\sigma_2
\intertext{and}
&\phi(\x)=
|\x|^2\sum_{i\ne 2} \frac{a_i}{(n+i)(i-2)}\sigma_i\big(\frac{\x}{|\x|}\big).
\end{align*}
\end{thm}

Let us explain how we are going to use Theorem \ref{KarpM} in the present paper: 
If $u$ is a solution to equation (\ref{main}) such that $u(0)=|\nabla u(0)|=0$,
if
\begin{align*}
&\lim_{j\to \infty}\frac{u(r_j\x)}{\sup_{B_{r_j}}|u|}=p
\intertext{for some $p\in \P$ and some sequence $r_j\to 0$, and if}
&\lim_{j\to \infty}\left( \frac{u(r_j\x)}{r_j^2}-\Pi(u,r_j)\right)=Z_p,
\end{align*}
then by $C^{1,\alpha}$-convergence we will have $Z_p(0)=|\nabla Z_p(0)|=0$
as well as $\Delta Z_p=-\chi_{\{p>0\}}$.
Also, by weak $W^{2,2}$-convergence we will have 
$$0=\int_{B_1} D^2 Z_p : D^2 h$$
for all $h\in \P$.
The latter is equivalent to  
$\Pi(Z_p,1)=0$.

By Theorem \ref{KarpM} with $\sigma = -\chi_{\{p>0\}}$ we can write
$$
Z_p(\x)=q(\x)\log|\x|+|\x|^2\phi(\x),
$$
where $q=\frac{a_2}{n+2}\sigma_2$. From here on we will assume that $n=3$ and $\x=(x,y,z)$.
It will also be convenient to parametrize the second order harmonic polynomials.
We will assume that 
$p=p_\delta  := (1/2+\delta)x^2+(1/2-\delta)y^2-z^2$. 
This can be done without loss of generality since there is always 
a rotation of the coordinate system such that $D^2 p$ is a diagonal matrix.
Rotating the coordinate system in that way, and if necessary
renaming $x,$ $y$ and $z$ we can always make sure that $p$ is 
up to a scaling factor
of the form
above or that $-p=(1/2+\delta)x^2+(1/2-\delta)y^2-z^2$. The latter
case can be handled similarly.
We would want to calculate $\sigma_2$. To that end we choose the polynomials
$3x^2-|\x|^2,$ $3y^2-|\x|^2$ and $3z^2-|\x|^2$
% , $xy$, $xz$ and $yz$ 
spanning
the axisymmetric second order harmonic polynomials in $\mathbb{R}^3$. 
That choice is somewhat arbitrary, but we contend that choosing different
polynomials would not facilitate substantially anything that follows.
It follows that
$$ 
\sigma_2=C\big(
(3A_x(\delta)-A(\delta))x^2+(3A_y(\delta)-A(\delta))y^2+
(3A_z(\delta)-A(\delta))z^2 \big),
$$
where
$C$ has been chosen such that
$\Vert \sigma \Vert_{L^2(\partial B_1)}=1$.

Using spherical coordinates
$x=r\sin(\theta)\cos(\phi),$ $y=r\sin(\theta)\sin(\phi)$ and
$z=r\cos(\theta)$, the characteristic function
$\chi_{\{p_\delta>0\}}=\chi_{\{ \theta > \textrm{arccot}(\sqrt{1/2+\delta\cos(2\phi)})\}}$.
The coefficients satisfy
\begin{align}
&A_x(\delta)=\int_{\partial
B_1}-\chi_{\{p_\delta>0\}}x^2=-8\int_{0}^{\pi/2}\int_{\textrm{arccot}(\sqrt{1/2+\delta
\cos(2\phi)})}^{\pi/2}\sin^3(\theta)\cos^2(\phi)d\theta d\phi,\non\\
&A_y(\delta)=\int_{\partial
B_1}-\chi_{\{p_\delta>0\}}y^2=-8\int_{0}^{\pi/2}\int_{\textrm{arccot}(\sqrt{1/2+\delta
\cos(2\phi)})}^{\pi/2}\sin^3(\theta)\sin^2(\phi)d\theta d\phi,\non\\
&A_z(\delta)=\int_{\partial
B_1}-\chi_{\{p_\delta>0\}}z^2=-8\int_{0}^{\pi/2}\int_{\textrm{arccot}(\sqrt{1/2+\delta
\cos(2\phi)})}^{\pi/2}\sin(\theta)\cos^2(\theta)d\theta d\phi\non\\
\intertext{and}
&A(\delta)=\int_{\partial
B_1}-\chi_{\{p_\delta>0\}}=-8\int_{0}^{\pi/2}\int_{\textrm{arccot}(\sqrt{1/2+\delta
\cos(2\phi)})}^{\pi/2}\sin(\theta)d\theta d\phi.\non
\end{align}

Next we notice that with
$${\bf K_0}=
\frac{2\log 2}{5 \Vert 3x^2-1\Vert^2_{L^2(\partial B_1)}},$$
\begin{equation}\label{piconst}
\Pi(Z_{p_\delta},1/2)=-{\bf K_0}\big(
(3A_x(\delta)-A(\delta))x^2+(3A_y(\delta)-A(\delta))y^2+
(3A_z(\delta)-A(\delta))z^2 \big),
\end{equation}
since $\Pi(\sigma_i)=0$ for $i\ne 2$.
Calculating $A_x,$ $A_y$, $A_z$
and $A$ we may estimate the rotation of $\Pi(u,r)$ as follows:
If $u(r\x )/r^2-\Pi(u,r)\approx Z_{p_\delta}$ then
$$
\Pi(u,r)-\Pi(u,r/2)\approx -\Pi(Z_{p_\delta},1/2).
$$
Later on this will be our main tool to analyze singular points. 

For convenience we will later also use the alternative representation
$p_\delta=(1-\delta)x^2+\delta y^2-z^2$ leading to the coefficients
\begin{align}
&B_x(\delta)=
-8\int_0^{\pi/2}\int_{\textrm{arccot}(\sqrt{(1-\delta)\cos^2(\phi)+\delta
\sin^2(\phi))})}^{\pi/2}
\sin^3(\theta)\cos(\phi)^2d\theta d\phi,\non\\
&B_y(\delta)=
-8\int_0^{\pi/2}\int_{\textrm{arccot}(\sqrt{(1-\delta)\cos^2(\phi)+\delta
\sin^2(\phi))})}^{\pi/2}
\sin^3(\theta)\sin(\phi)^2d\theta d\phi\non\\
\intertext{and}
&B(\delta)=
-8\int_0^{\pi/2}\int_{\textrm{arccot}(\sqrt{(1-\delta)\cos^2(\phi)+\delta
\sin^2(\phi))})}^{\pi/2}
\sin(\theta)d\theta d\phi.\non
\end{align}
That is, $B_x(1/2-\delta)=A_x(\delta)$ etc.
It will be convenient to define the parameter $\delta$ for polynomials
and solutions:
\begin{definition}\label{delta}
For each $v\in W^{2,2}(B_1)$, let, if necessary after rotation,
\begin{align}
\Pi(v,1)/\sup_{B_1} |\Pi(v,1)|=&(1/2+\delta)x^2+(1/2-\delta)y^2-z^2\non\\
\textrm{or } &-(1/2+\delta)x^2-(1/2-\delta)y^2+z^2.\non
\end{align}
We note that $\delta$ is unique and
define $\delta^A(v) := \delta$.
\\
Moreover, let
\begin{align}\Pi(v,1)/|\sup_{B_1} \Pi(v,1)|=&(1-\tilde\delta)x^2+\tilde\delta y^2-z^2\non\\
\textrm{or } &-(1-\tilde\delta)x^2-\tilde\delta y^2+z^2.\non
\end{align}
We note that $\tilde \delta$ is unique and
define $\delta^B(v) := \tilde\delta$.
It is important to note that $\sup_{B_1}|\Pi(v,1)|/\sup_{B_1} |\Pi(v,1)|=1$
implies $\delta^B(v)=\tilde \delta\ge 0$.
\\
We will also use $\delta^A_r := \delta^A(u(r\cdot))$
and $\delta^B_r := \delta^B(u(r\cdot))$.
\end{definition}
In general we cannot calculate the integrals $A_x,A_y,\dots$ explicitly.
In some special cases however, when we have sufficient symmetry, we may even write down
explicit solutions to the equation $\Delta u=-\chi_{\{p>0\}}$.
Luckily and surprisingly, as seen in Section \ref{SecClass}, these special solutions
are the only solutions appearing as limits of
$$\frac{u(r\x+\x^0)}{r^2}-\Pi(u,r,\x^0).$$

{\bf 1. $2$-dimensional Solutions} (cf. \cite[Lemma 4.4]{cmp}):

Define $v:(0,+\infty)\times [0,+\infty)\to \mathbb{R}$ by
$$
v(x,z):=
-4xz\log(x^2+z^2)+2(x^2-z^2)\left( \frac{\pi}{2}-2\arctan\left(\frac{z}{x}\right)
\right)-\pi(x^2+z^2).
$$
Moreover, let
$$
w(x,z):=\left\{ \begin{array}{ll}
v(x,z),& xz\ge 0, x\ne 0,\\
-v(-x,z),& x<0, z\ge 0,\\
-v(x,z),& x>0, z\le 0,
\end{array}\right.$$
and define
$$
Z(x,z) := \frac{w(x,z)-\pi(x^2+z^2)+8xz}{8\pi}.
$$
It has been shown in \cite[Lemma 4.4]{cmp} that\\
(i) $\Delta Z=-\chi_{\{xz>0\}}$ in $\mathbb{R}^3$.\\
(ii) $Z(0)=|\nabla Z(0)|=0$.\\
(iii) $\lim_{|\x| \to \infty}{Z(\x)\over {|\x|^3}}=0$.\\
(iv) $\Pi(Z,1)=0$.\\
(v) $\Pi(Z,1/2)=\log(2)xz/\pi, \;\tau(Z,1/2)=\log(2)/(2\pi)$.\\
(vi) $Z$ is the unique function satisfying (i)-(iv).

$ $

{\bf 2. True 3D Solutions}:

Next we are going to calculate $Z_p$
for $p(\x)=(x^2+y^2)/2-z^2$.

Let us denote
\begin{align*}
&v_1=2p(\x)\log\big( x^2+y^2 \big)-4z^2
\intertext{and}
&v_2=-\frac{3z|\x|}{2}+\frac{1}{2}p(\x)
\log\Big(\frac{|\x|-z}{|\x|+z} \Big).
\end{align*}
Then $\Delta v_1=0$ and $\Delta v_2=0$ in 
$\mathbb{R}^3_+\setminus\{x^2+y^2=0\}$. Also notice that 
$\partial_z v_1(x,y,0)=0$.

Let 
$$
v(\x)=\frac{\sqrt{3}}{36}\Big( 4v_2(\x)-v_1(\x)\Big)
$$
in $\{p(\x)\le 0\}\cap \{z>0\}$, where the coefficients for $v_1$ and
$v_2$ are chosen such that the singularities cancel at $x=y=0$.
Moreover, let
\begin{align}
v(\x)=&-\frac{\sqrt{3}}{36}v_1(\x)\non\\
&+\frac{3+\sqrt{3}\log\big(2-\sqrt{3}\big)}{18}p(\x)-\frac{|\x|}{6}
\non\end{align}
in $\{p(\x)>0\}\cap \{z>0\}$.

Next we reflect $v$ at $\{z=0\}$ according to
$$
\tilde{Z}_1(\x)=\left\{
\begin{array}{ll}
v(x,y,z) & \textrm{ if } z\ge 0, \\
v(x,y,-z) & \textrm{ if } z<0.
\end{array}\right.
$$ 
Last, we define $Z_1=\tilde{Z}_1-\Pi(\tilde{Z}_1,1)$
and $Z_2=-Z_1$.

We have thus established the following lemma:
\begin{lem}\label{zcalculations}
Let $Z$, $Z_1$ and $Z_2$ be as above. Then, with $p(\x)=(x^2+y^2)/2-z^2$,
\begin{enumerate}
\item $\Delta Z=-\chi_{\{xz>0\}}$, $Z(0)=|\nabla Z(0)|=0$,
$\lim_{\x\to \infty}|Z(\x)|/|\x|^3=0$ and
$\Pi(Z,1)=0, \Pi(Z,1/2)=(\log(2)/\pi)xz$.
\item $\Delta Z_1=-\chi_{\{p(\x)>0\}}$, $Z_1(0)=|\nabla Z_1(0)|=0$, \\
$\lim_{\x\to \infty}|Z_1(\x)|/|\x|^3=0$\\ and
$\Pi(Z_1,1)=0, \Pi(Z_1,1/2)=\log(2)(\sqrt{3}/9)p(\x)$.
\item $\Delta Z_2=-\chi_{\{p(\x)<0\}}$, $Z_2(0)=|\nabla Z_2(0)|=0$, \\
$\lim_{\x\to \infty}|Z_2(\x)|/|\x|^3=0$\\ and
$\Pi(Z_2,1)=0, \Pi(Z_2,1/2)=-\log(2)(\sqrt{3}/9)p(\x)$.
\end{enumerate}
\end{lem}
\proof The proof follows from simple calculation.\qed
\begin{remark}
The fact that $\Pi(Z_p,1/2)$ is a multiple of the polynomial $p$ in the above
three cases, natural though it may be, will be of paramount importance
in later chapters when it comes to the question of unique tangent cones.
\end{remark}
The following two collections of properties of the A's and B's
visualized in Figure \ref{fig:figure1}-\ref{fig:figure5} are of central importance in our paper and will be proved in the Appendix together with Lemma \ref{etalem} below.
\begin{figure}[ht]
\begin{minipage}[t]{0.45\linewidth}
\centering
\includegraphics[width=6cm]{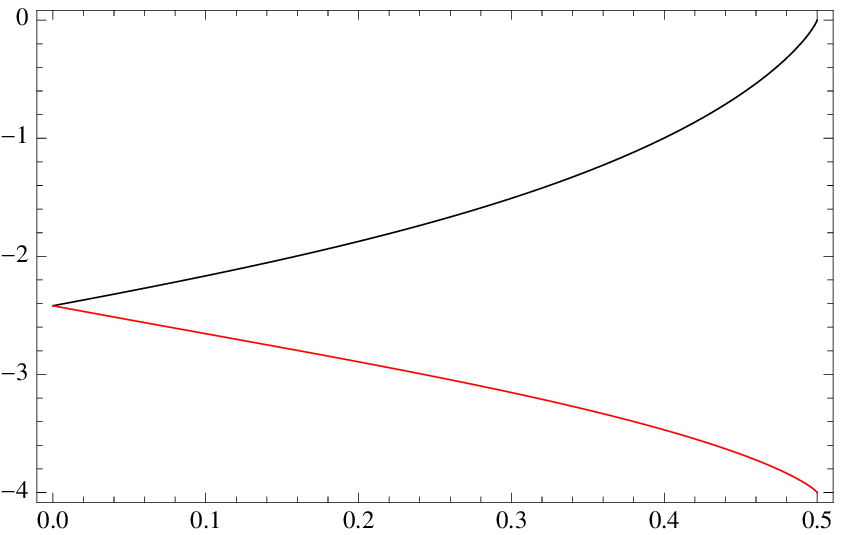}
\caption{$3A_y(\delta)-A(\delta)$ (upper graph) versus $3A_x(\delta)-A(\delta)$}
\label{fig:figure1}
\end{minipage}
\hspace{0.5cm}
\begin{minipage}[t]{0.45\linewidth}
\centering
\includegraphics[width=6cm]{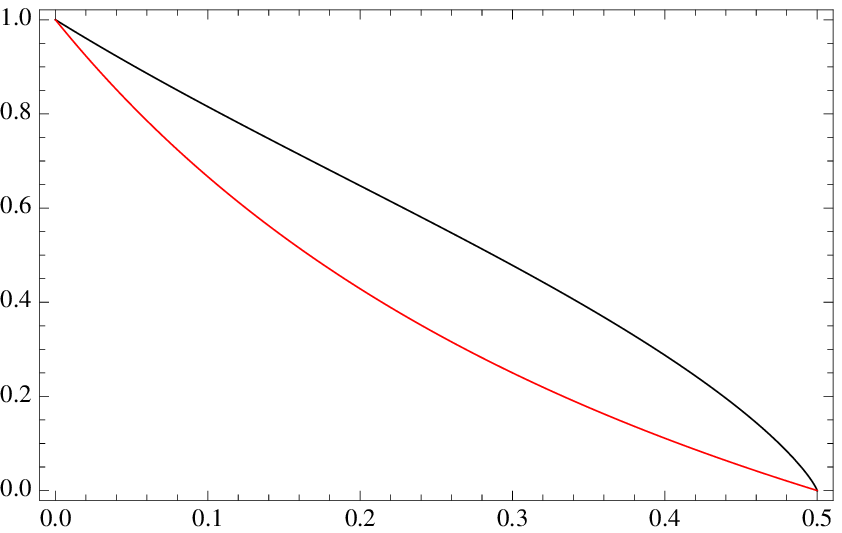}
\caption{$\frac{3A_y(\delta)-A(\delta)}{3A_x(\delta)-A(\delta)}$ (upper graph) versus
$\frac{1-2\delta}{1+2\delta}$}
\label{fig:figure2}
\end{minipage}
\end{figure}
\begin{thm}\label{athm}
For $\delta\in (0,1/2)$:\\
(i) $(3A_y(\delta)-A(\delta))'>0$.\\
(ii) $(3A_x(\delta)-A(\delta))'<0$.\\
(iii) $3A_x(0)-A(0)=3A_y(0)-A(0)<0$.\\
(iv) $3A_x(\delta)-A(\delta)<0$.\\
(v) $3A_y(1/2)-A(1/2)=0$.\\
(vi) $3(A_y''-A_x'')+2\delta (3A_y''+3A_x''-2A'')+2(3A_y'+3A_x'-2A')>0$.\\
(vii) $$
\frac{3A_y(\delta)-A(\delta)}{3A_x(\delta)-A(\delta)}>
\frac{1-2\delta}{1+2\delta}.$$
(viii) $3A_x'(0)-A'(0)= -\frac{4\pi}{3\sqrt{3}}$.\\
(ix) $3A_y'(0)-A'(0)=\frac{4\pi}{3\sqrt{3}}.$
\end{thm}
\begin{lem}\label{etalem}
For the positive universal constant
$\eta_0$ defined in (\ref{eta0}) and every
$\delta\in (0,1/2)$,
$$
\sup_{B_1}|Cp_\delta +\Pi(Z_{p_\delta},1/2)|\ge C+\eta_0,
$$
for every sufficiently large constant $C>0$.
\end{lem}
\begin{figure}[ht]
\begin{minipage}[t]{0.45\linewidth}
\centering
\includegraphics[width=6cm]{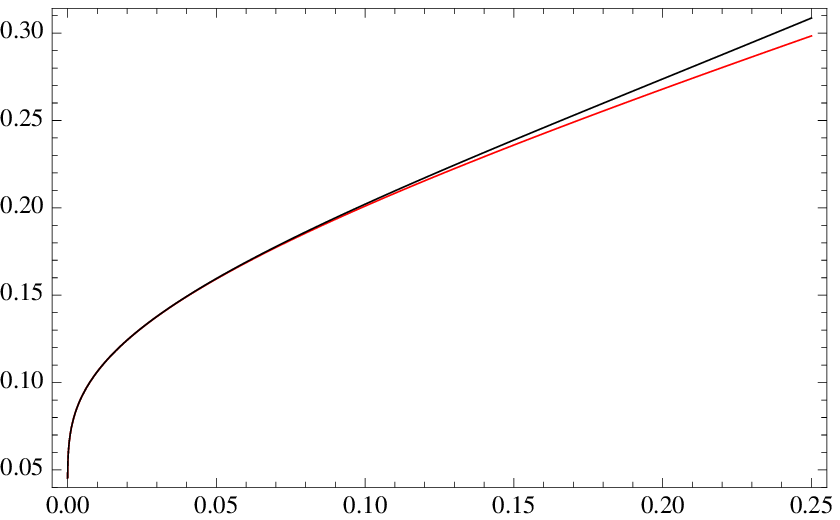}
\caption{$1/C_0$ vs. $1/C_y$}
\label{fig:figure3}
\end{minipage}
\hspace{0.5cm}
\begin{minipage}[t]{0.45\linewidth}
\centering
\includegraphics[width=6cm]{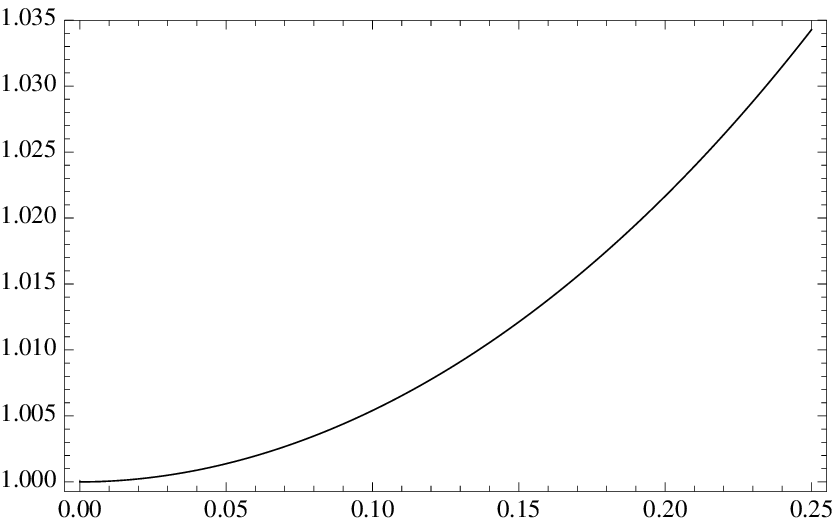}
\caption{$C_y/C_0$}
\label{fig:figure4}
\end{minipage}
\hspace{0.5cm}
\begin{minipage}[b]{0.45\linewidth}
\centering
\includegraphics[width=6cm]{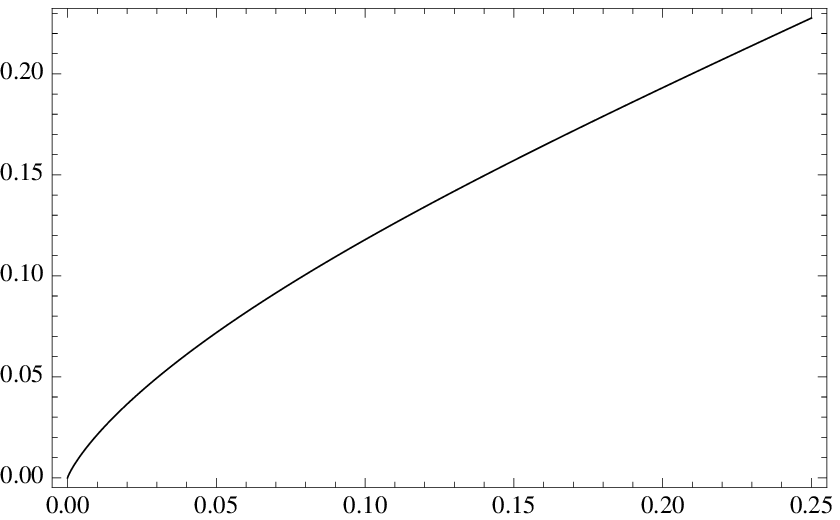}
\caption{$C_x(\delta)/\delta$}
\label{fig:figure5}
\end{minipage}
\end{figure}
\begin{thm}\label{bthm}
For small $\delta>0$:\\
(i) $
B_y(\delta)=B_y(0)-C_y(\delta)\delta,
$
where $0<k_1 \le C_y(\delta)/|\log\delta|\le K_2<+\infty$. \\
(ii) $
B(\delta)=B(0)-C_0(\delta)\delta,
$
where $\lim_{\delta\to 0}C_0(\delta)=+\infty$.\\
(iii) $
B_x(\delta)=B_x(0)+o(\delta)$ as $\delta\to 0$.
\\
(iv) $
\lim_{\delta\to 0}\frac{C_y(\delta)}{C_0(\delta)}=1.
$
\end{thm}
\section{Growth of the Solution}\label{growth}
Since $Z, Z_1, Z_2$ have growth $|\x|^2|\log(|\x|)|$ away from the origin
and we expect $u$ (up to some harmonic part)
to be close to $Z_p$, $u$ should share the same growth. We will
prove
this in the next lemma.

Before we state the lemma let us point out a simple fact that will
be used frequently in what follows. By Proposition \ref{boundfordiff}
we know that if the origin is a singular point then $u(r\x )/r^2-\Pi(u,r)$
is uniformly bounded by a constant $C_0$ depending on $n$ and $\|u\|_{L^\infty(B_1)}$.
This implies that when $u(r\x )/r^2$ is large, say $
\max(\sup_{B_1}|\Pi(u,r)|,\sup_{B_1}|u(r\x )/r^2|)\ge 2C_0$, then
\begin{equation}\label{comp}
\frac{1}{2}\sup_{B_1}|\Pi(u,r)|\le \sup_{B_1}|u(r\x )/r^2|\le 2\sup_{B_1}|\Pi(u,r)|.
\end{equation}
So controlling the size of $\Pi(u,r)$ is equivalent to controlling the size of $u(r\x)/r^2$
at singular points.

In the two-dimensional case the following lemma has been proved in \cite[Lemma 5.5]{cmp}.
\begin{lem}\label{superquad} Let $n=3$ and let $u$ be a solution to (\ref{main}) in $B_1$
such that $\sup_{B_1}|u|\le M$ and
$u(0)=|\nabla u(0)|=0$.
Then there exist $\rho_0>0$ and $r_0>0$ such that if 
\begin{equation}\label{delta0}
\sup_{B_1}|\Pi(u,r)|\ge \frac{1}{\rho_0}
\end{equation}
for an $r\le r_0$ then
$$
\sup_{B_1}|\Pi(u,r/2)|>\sup_{B_1}|\Pi(u,r)|+\eta_0/2,
$$
where $\eta_0$ is the positive constant in Lemma \ref{etalem}
(see (\ref{eta0}) in the Appendix).
\end{lem}
\proof
If the Lemma is not true, then there exists a sequence $u^j$ of solutions to
(\ref{main}) and  $r_j\to 0$ such that
$$
\sup_{B_1}|\Pi(u^j,r_j)|\ge j
\textrm{ and }\sup_{B_1}|\Pi(u^j,r_j/2)|\le \sup_{B_1}|\Pi(u^j,r_j)|+\eta_0/2.$$

Using Proposition \ref{boundfordiff}, and passing if necessary to a subsequence,
$$
v^j(\x):=\frac{u^j(r_j \x)}{r_j^2}-\Pi(u^j,r_j)\to v \textrm{ in }
C^{1,\alpha}_{\rm loc}(\mathbb{R}^3)\cap W^{2,p}_{\rm loc}(\mathbb{R}^3).
$$
We also have $\Pi(u^j,r_j/2)=\Pi(u^j,r_j)+\Pi(v^j,1/2)$.
The limit $v$ satisfies
$\Delta v= -\chi_{\{p>0\}}$, where 
---passing if necessary to another subsequence---
 $$p=\lim_{j\to \infty} p_j \textrm{ for } p_j := \Pi(u^j,r_j)/\|\Pi(u^j,r_j)\|_{L^\infty (B_1)}.$$
It follows that $\Pi(v,1/2)=\Pi(Z_p,1/2)$ where $Z_p$ is the unique
solution to
$$
\begin{array}{ll}
\Delta Z_p=-\chi_{\{p>0\}} & \textrm{ in } \mathbb{R}^3, \\
Z_p(0)=|\nabla Z_p(0)|=\Pi(Z_p,1)=0 & \textrm{ and } \\
\lim_{|\x|\to \infty}\frac{|Z_p(\x)|}{|\x|^3}=0. &
\end{array}
$$
Consequently
$\lim_{j\to \infty}(\Pi(u^j,r_j/2)-\Pi(u^j,r_j))=\Pi(Z_p,1/2)$, and
\begin{equation}\label{jui}
\sup_{B_1}|\Pi(u^j,r_j/2)|=\sup_{B_1}|\Pi(u^j,r_j)+\Pi(Z_{p_j},1/2)|+o(1)
\textrm{ as }j\to\infty.
\end{equation}
Finally we apply Lemma \ref{etalem}
and obtain the statement of the lemma. 
\qed
\begin{remark}\label{hdimlem}
Lemma \ref{superquad} extends to dimension $n>3$ provided that
for some $\epsilon$ depending only on $n$ and $M$
$$\sup_{B_1} \left|\frac{\Pi(u,r)}{\sup_{B_1} |\Pi(u,r)|}-p(Q\cdot)\right|
\le \epsilon$$
for a three dimensional polynomial $p\in \P$
and a rotation $Q$ in $\mathbb{R}^n$.
\end{remark}
In the two-dimensional case the following lemma has been proved in \cite[Corollary 5.6]{cmp}.
\begin{cor}\label{loggrowth}
Let $n=3$ and let $u$ be a solution to (\ref{main}) in $B_1$
such that $\sup_{B_1}|u|\le M$ and
$u(0)=|\nabla u(0)|=0$.
Then there exist $\rho_0>0$ and $r_0>0$ such that if
\begin{align*}
&\sup_{B_1}|\Pi(u,r)|\ge \frac{1}{\rho_0}
\intertext{for an $r\le r_0$ then}
&\sup_{B_1}|\Pi(u,2^{-j}r)|\ge \sup_{B_1}|\Pi(u,r)|+j \eta_0/2
\intertext{and}
&\sup_{B_s}|u|\ge \frac{1}{16}\bigg(\Big( \frac{s}{r}\Big)^2\sup_{B_r}|u|+\eta_0 s^2\log(r/s)\bigg)\textrm{ for } 0<s<r,
\end{align*}
where $\eta_0$ is the positive constant in Lemma \ref{etalem}.

Furthermore, there exists a constant $\kappa=\kappa(M,n)$ such that
\begin{align*}
&\sup_{B_1}|\Pi(u,2^{-j}r)|\le \sup_{B_1}|\Pi(u,r)| + \kappa j,\\
&
\sup_{B_1}|\Pi(u,s)|\le 2\sup_{B_1}|\Pi(u,r)|+\kappa j\textrm{ for }s\in (2^{-(j+1)}r, 2^{-j}r],
\intertext{and}
&\sup_{B_s}|u|\le 16\bigg(\Big( \frac{s}{r}\Big)^2\sup_{B_r}|u|+\kappa s^2 \log(r/s)\bigg).
\end{align*}
\end{cor}
\proof
Lemma \ref{superquad} applies, so that
$$
\sup_{B_{1}}|\Pi(u,2^{-1}r)|\ge \sup_{B_1}|\Pi(u,r)|+\eta_0/2\ge \frac{1}{\rho_0}.
$$
It follows that Lemma \ref{superquad} applies again with $2^{-1}r$.
Thus we may iterate
and deduce that
$$
\sup_{B_1}|\Pi(u,2^{-j}r)|\ge \sup_{B_1}|\Pi(u,r)|+j\eta_0/2,
$$
which together with (\ref{comp}) proves the first part of the Corollary.

We also notice that by Proposition \ref{boundfordiff} we have
$$
\sup_{B_1}|\Pi(u,r/2)|\le \sup_{B_1}|\Pi(u(r\x )/r^2,1/2)-\Pi(u,r)|+\sup_{B_1}|\Pi(u,r)|\le
\kappa +\sup_{B_1}|\Pi(u,r)|.
$$
Arguing as above
we get
$$
\sup_{B_1}|\Pi(u, 2^{-j}r)|\le \sup_{B_1}|\Pi(u,r)|+j\kappa,
$$
which together with (\ref{comp}) proves the second part of the Corollary.\qed
\section{Existence of a True Three-dimensional Singularity}\label{truethree}
\begin{cor}\label{exist}
There exists a solution $u$ of (\ref{main}) in $B_1\subset \mathbb{R}^3$ such that
$$
\lim_{r\to 0}\frac{u(r\x )}{\sup_{B_r}|u|}
=  \frac{x^2+y^2}{2}-z^2 
$$
or
$$
\lim_{r\to 0}\frac{u(r\x )}{\sup_{B_r}|u|}
=  z^2-\frac{x^2+y^2}{2}. 
$$
\end{cor}
\proof The proof is similar to
that of
\cite{anderssonweiss}, so we will only give a sketch.
We define the operator
$T=T_{\epsilon}:C^{\alpha}(B_1^+)\to C^{\alpha}(B_1^+)$ by
$$
\begin{array}{ll}
\Delta T(u)=-f_{\epsilon}(u-u(0)) & \textrm{ in } B_1^+, \\
T(u)= M\big(\frac{x^2+y^2}{2}-z^2\Big) & \textrm{ on } \partial B_1\cap \{ z>0\}, \textrm{ and}\\
\frac{\partial T(u)}{\partial z}=0 & \textrm{ on } \{z=0\}\cap B_1.
\end{array}
$$
Moreover we impose that $T(u)$ has cylindrical symmetry, that is
$T(u)(x,y,z)=g(x^2+y^2,z)$
for some function $g$. The function $f_\epsilon(t)$ is a smooth
approximation of $\chi_{\{t>0\}}$ and $M$ is some large constant.

By Schauder's fixed point theorem there exists an $u_\epsilon$
such that $T_\epsilon(u_\epsilon)=u_\epsilon$. We may pass to the limit
$\lim_{\epsilon \to 0}u_\epsilon=\tilde{u}$. Defining
$u(\x)=\tilde{u}(\x)-\tilde{u}(0)$ for $z>0$ and
$u(\x)=\tilde{u}(x,y,-z)-\tilde{u}(0)$ for $z<0$, we see that $u$ solves
(\ref{main}). 
From the boundary condition we infer as in \cite{anderssonweiss}
that $\Phi^u_{0}(r)\le -M$, which implies 
that $\sup_{B_1}|\Pi(u,r)|$ is also large.
From Corollary \ref{loggrowth} we conclude therefore that
$$
\sup_{B_s}|u|\ge \frac{1}{16}\bigg(\Big( \frac{s}{r}\Big)^2(\sup_{B_r}|u|+\eta_0s^2 \log(r/s)\bigg).
$$
But then 
each limit of 
$$\frac{u(r\x )}{\sup_{B_r}|u|}$$
as $r\to 0$ must be a polynomial $p\in \P$.
Naturally, $p$
will have the same cylindrical symmetry as $u$. Therefore
$p=(x^2+y^2)/2-z^2$ or $p=z^2-(x^2+y^2)/2$. 
\\
Last suppose towards a contradiction that there are two subsequences
such that one converges to $(x^2+y^2)/2-z^2$
and the other to $z^2-(x^2+y^2)/2$.
By a continuity argument we obtain in this case a third subsequence
and a limit that is neither $(x^2+y^2)/2-z^2$
nor $z^2-(x^2+y^2)/2$, a contradiction.
\qed
\section{Estimating $u-\Pi(u)-Z_{\Pi(u)}$}\label{prelanal}

The following Lemma is a direct consequence of Corollary 4.1 in 
\cite{Ganzburg}.

\begin{lem}\label{Cuyt}
Let $p$ be a second order polynomial in $\mathbb{R}^n$ and
$\|p\|_{L^{\infty}(Q_1)}=1$. Then
$$
\Big(\frac{\L^n(\{|p|\le \epsilon\})}{|\log(\L^n(\{|p|\le
\epsilon\}))|^{n-1}} \Big)^2\le C(n)\epsilon
\textrm{ for every } \epsilon \in (0,1).
$$
In particular,
$$
\L^{n}(\{|p|\le \epsilon\})\le C(n,\alpha)\epsilon^{\alpha}\textrm{ for every } \epsilon \in (0,1)
$$
and all $\alpha<1/2$.
\end{lem}
The following Lemma is related to the two-dimensional result \cite[Lemma 6.1]{cmp}.
\begin{lem}\label{gest}
Let $u$ solve (\ref{main}) in $B_1\subset \R^n$ such that $\sup_{B_1}|u|\le M$ and
$u(0)=|\nabla u(0)|=0$, and for some $\rho\le \rho_0$ and $r\le r_0$ let
$$
\sup_{B_1}|\Pi(u,r)|\ge \frac{1}{\rho}.
$$
Furthermore let $g_r$ be the solution of
$$
\begin{array}{ll}
\Delta g_r=\chi_{\{\Pi(u,r)>0\}}-\chi_{\{u(r\cdot)>0\}} & \textrm{ in }B_1,\\
g_r=0 & \textrm{ on } \partial B_1.
\end{array}
$$ Then for each $\alpha<1/4$,
\begin{enumerate}
\item $\Vert D^2 g_r\Vert_{L^2(B_1)} \le C(M,n,\alpha) \big| \sup_{B_1}|\Pi(u,r)|\big|^{-\alpha}$.
\item $\max\Big(\sup_{B_1} |\Pi(g_r,1)|,\sup_{B_1} |\Pi(g_r,1/2)|\Big)\le C(M,n,\alpha) \big|\sup_{B_1}|\Pi(u,r)|\big|^{-\alpha}.$
\end{enumerate}
\end{lem}
\proof
Let $p=\Pi(u,r)$.
We know that $\Delta g_r=1$ when $p>0$ and $u(r\x)\le 0$,
and that $\Delta g_r=-1$ when $p\le 0$ and $u(r\x)>0$; in all other cases it is $0$. By
Proposition \ref{boundfordiff} we also have that 
$$
\Big| \frac{u(r\x)}{r^2}-p\Big|\le C_0.
$$
Combining those properties we obtain that $\Delta g_r=0$ outside the
set $\{|p|\le C_0\}$.
From Lemma \ref{Cuyt} it follows that
$$\|\Delta g_r\|_{L^2(B_1)}\le \left(\L^n(\{|p|\le C_0\})\right)^{1\over 2}\le C(M,n)|\sup_{B_1}|p||^{-\alpha} \textrm{ for each }\alpha<1/4.$$
Standard $L^2$-theory (see for example \cite{Stein})
thus implies (i).

Rotating and setting $q:= \Pi(g_r,t) = \sum_{j=1}^n a_j x_j^2$,
where $t=1$ or $t=1/2$,
we obtain
\begin{align*}
&\Vert D^2 q\Vert_{L^2(B_1)} \le C_1 \Vert D^2 g_r\Vert_{L^2(B_1)}
\le C_2 \big|\sup_{B_1}|\Pi(u,r)|\big|^{-\alpha}
\intertext{and}
&|a_j| \le C_3  \big|\sup_{B_1}|\Pi(u,r)|\big|^{-\alpha}
\end{align*}
for every $1\le j \le n$, proving (ii).\qed

The following Lemma is related to the two-dimensional result \cite[Lemma 4.3]{cmp}.
\begin{cor}\label{corgest}
Let $u$ solve (\ref{main}) in $B_1\subset \R^n$ and
assume that $\sup_{B_1}|u|\le M$,
$u(0)=|\nabla u(0)|=0$, and that for some $\rho\le \rho_0$ and $r\le r_0$,
$$
\sup_{B_1}|\Pi(u,r)|\ge \frac{1}{\rho}.
$$
Then
$$
\sup_{B_1}|\Pi(u,r/2)-\Pi(u,r)-\Pi(Z_{\Pi(u,r)},1/2)|\le C(M,n,\alpha)(\sup_{B_1}|\Pi(u,r)|)^{-\alpha}
$$
for each $\alpha<1/4$.
\end{cor}
\proof  For each $r$ write
$$
u(r\x)/r^2=\Pi(u,r)+Z_{\Pi(u,r)}+\tilde{g}_r+\tilde{h}_r
$$
where
$
\Delta \tilde{g}_r= \Delta (u(r\x)/r^2)-\Delta Z_{\Pi(u,r)},
$
$
\Delta \tilde{h}_r=0
$
and
$\tilde{g}_r(0)=|\nabla
\tilde{g}_r(0)|=\tilde{h}_r(0)=|\nabla\tilde{h}_r(0)|=|\Pi(\tilde{g}_r,1)|=|\Pi(\tilde{h}_r,1)|=0$.
Next denote by $g_r$ the solution to
$$
\begin{array}{ll}
\Delta g_r =\Delta \tilde{g}_r & \textrm{ in } B_1, \\
g_r=0 & \textrm{ on } \partial B_1.
\end{array}
$$
Then $\tilde{g}_r+\tilde{h}_r=g_r+h_r$ for some harmonic function $h_r$ in $B_1$. From Lemma
\ref{gest} it follows that
$$
\sup_{B_1}|\Pi(g_r,1/2)|\le C_1(M,n,\alpha)(\sup_{B_1}|\Pi(u,r)|)^{-\alpha},
$$
and from Lemma \ref{proj} we infer that
\begin{align*}
\sup_{B_1}|\Pi(h_r,1/2)|&=\sup_{B_1}|\Pi(h_r,1)|=\sup_{B_1}|\Pi(\tilde{g}_r,1)-\Pi(g_r,1)|\\&\le
C_1(M,n,\alpha)(\sup_{B_1}|\Pi(u,r)|)^{-\alpha}.
\end{align*}
Thus in $B_1$
\begin{align*}
&\Pi(u,r/2)=\Pi(u,r)+\Pi(Z_{\Pi(u,r)},1/2)+\Pi(g_r,1/2)+\Pi(h_r,1/2)
\intertext{and}
&|\Pi(u,r/2)-\Pi(u,r)-\Pi(Z_{\Pi(u,r)},1/2)|\le C(M,n,\alpha)(\sup_{B_1}|\Pi(u,r)|)^{-\alpha}.\qed
\end{align*}
\section{Classification of Blow-up Limits in $\mathbb{R}^3$ \\ --- An Unexpected Symmetrization Effect}\label{SecClass}

In this section we will show that if
$\lim_{j\to \infty}\frac{u(r_j\x)}{\sup_{B_{r_j}}|u|}=p$, where $p$ is a harmonic
polynomial,
then $p=2 xz$, $p=(x^2+y^2)/2-z^2$ or $p=z^2-(x^2+y^2)/2$ up to a
rotation.
\begin{thm}\label{symm}
Let $n=3$, let $\Delta u=-\chi_{\{u>0\}}$ in $B_1$ and assume that $u(0)=|\nabla u(0)|=0$ and that the
monotonicity energy satisfies $\lim_{r\to 0}\Phi^u_{\x^0}(r)=-\infty$. Then
each limit of 
$$
\frac{u(r\x)}{r^2}-\Pi(u(\x),r),
$$
as $r\to 0$, is contained in
$$
\{Z(Q\x):\; Q\in \mathcal{Q}\}\cup \{Z_1(Q\x):\; Q\in \mathcal{Q}\}\cup \{Z_2(Q\x):\;
Q\in \mathcal{Q}\};
$$
here $\mathcal{Q}$ is the set of all rotations of $\mathbb{R}^3$.
\end{thm}
\proof Suppose towards a contradiction that the statement is not true.
By Proposition \ref{fixedcenter} (i) 
there exists a solution $u$ and a sequence $r_j\to 0$ such that after rotation
\begin{equation}\label{symass}
\lim_{j\to \infty}\Big| \frac{u(r_j\x)}{\sup_{B_{r_j}}|u|}-p_{\delta_0}\big|=0
\end{equation}
for some $\delta_0\in (0,1/2)$ 
and $p_{\delta_0}=(1/2+\delta_0)x^2+(1/2-\delta_0)y^2-z^2$ and
$\delta_0\in (0,1/2)$ or
$p_{\delta_0}=z^2-(1/2+\delta_0)x^2-(1/2-\delta_0)y^2$.
We may assume that $p_{\delta_0}=(1/2+\delta_0)x^2+(1/2-\delta_0)y^2-z^2$.
Furthermore, from Proposition \ref{fixedcenter} (i) and Proposition
\ref{boundfordiff}, $\lim_{r\to 0}\sup_{B_1}|u(r\x)/r^2|=\infty$.

We are going to prove a decay estimate
for $\delta^A_r=\delta^A(u(r\cdot))$ in $r$ which will lead to a contradiction
to (\ref{symass}).

By Theorem \ref{athm} (vii),
$$
\frac{3A_y(\delta)-A(\delta)}{3A_x(\delta)-A(\delta)}>
\frac{1-2\delta}{1+2\delta}
\textrm{ for } \delta\in (0,1/2).
$$
Thus
\begin{equation}\label{gammagain10}
\kappa(\delta) := (1+2\delta)
\frac{3A_y(\delta)-A(\delta)}{3A_x(\delta)-A(\delta)}-
(1-2\delta) \ge \omega(\delta) > 0,
\end{equation}
where $\omega$ is a continuous function on $[0,1/2]$.

By Corollary \ref{corgest}, 
using Corollary \ref{loggrowth} to estimate $$(\sup_{B_1}|\Pi(u,r)|)^{-\alpha}=O(|\log(r)|^{-\alpha}),$$ we obtain
for every
$\alpha<1/4$ that in $B_1$, up to a rotation depending on $r$,
\begin{equation}\label{bigO}
|\Pi(u,r/2)-\tau_r p_{\delta^A_r}-\Pi(Z_{\delta^A_r},1/2)|\le O(|\log(r)|^{-\alpha});
\end{equation}
from here on, $Z_{\delta}$ is the unique solution to
$$
\begin{array}{ll}
\Delta Z_{\delta}=-\chi_{\{p_{\delta}>0\}} & \textrm{ in } \mathbb{R}^3 \\
Z_{\delta}(0)=|\nabla Z_{\delta}(0)|=\Pi(Z_{\delta},1)=\lim_{|\x|\to\infty}Z_{\delta}(\x)/|\x|^3=0. & 
\end{array}
$$
In particular, for the ${\bf K_0}$ defined in (\ref{piconst}),
\begin{align*}
\Pi(Z_{\delta},1/2)&=-
{\bf K_0}\big[
(3A_{x}(\delta)-A(\delta))x^2+(3A_y(\delta)-A(\delta))y^2+(3A_z(\delta)-A(\delta))z^2\big]\\
&=-{\bf K_0} K_1 [(1+2\delta)x^2 + K_2 (1+2\delta)y^2 + K_3 (1+2\delta) z^2],
\end{align*}
where --- using the fact that $\Pi(Z_{\delta},1/2)$ is harmonic ---
\begin{align*}
&K_1= \frac{3A_x(\delta)-A(\delta)}{(1+2\delta_0)},\\
&K_2 = \frac{3A_y(\delta)-A(\delta)}{3A_x(\delta)-A(\delta)}
= \frac{1-2\delta+\kappa(\delta)}{1+2\delta}\\
\intertext{and}
&K_3 = \frac{3A_z(\delta)-A(\delta)}{3A_x(\delta)-A(\delta)}
= -(1+K_2) = -\frac{2+\kappa(\delta)}{1+2\delta}.
\end{align*}
It follows that
\begin{align*}\Pi(Z_{\delta},1/2)
&= -
{\bf K_0}K_1 [(1+2\delta)x^2 + (1-2\delta)y^2 + \kappa(\delta) y^2 - 2z^2 -\kappa(\delta) z^2]\\
&= c_\delta \big( 2p_{\delta}+\kappa(\delta) (y^2-z^2)\big)
\end{align*}
for $c_\delta =-{\bf K_0}(3A_x(\delta)-A(\delta))/ (1+2\delta)\ge \bar c>0, \delta\in [0,1/2]$
(see Theorem \ref{athm}).

Invoking (\ref{bigO}), this implies that in $B_1$, up to a rotation depending on $r$,
\begin{equation}\label{firststep}
|\Pi(u,r/2)-\tau_r p_{\delta^A_r}-c_{\delta^A_r}\big( 2p_{\delta^A_r}+\kappa(\delta^A_r) (y^2-z^2)\big)|=O(|\log(r)|^{-\alpha}).
\end{equation}

The fact that $c_\delta\ge \bar c>0$
as well as the estimate $\kappa\ge \omega>0$
consequently prove together with Corollary \ref{loggrowth} that, rotating slightly,
\begin{equation}\label{it1}
\delta^A_{r/2}\le \delta^A_{r} - c_1 
\frac{\omega(\delta^A_{r})}{|\log r|}+C_2 |\log(r)|^{-1-\alpha}.
\end{equation}
Note that estimate (\ref{it1}) is independent of rotations.
As long as $\omega(\delta^A_{2^{-k}r_0})\ge  |\log(2^{-k}r_0)|^{-\alpha/2}$,
induction of estimate (\ref{it1}) in $k$ yields a logarithmic decay
of $\delta^A_r$ in $r$.
It follows that $\delta^A_r\to 0$ as $r\to 0$,
contradicting
the assumption
$\lim_{j\to\infty}\delta^A_{r_j}=\delta_0>0$.\qed
\section{Unique Tangent Cones at True Three Dimensional Singularities}\label{uniquesingpoints}

From the previous section we may infer by a continuity argument that 
in three dimensions, assuming $u(0)=|\nabla u|=0$ as well as
$\lim_{r\to 0}\Phi^u_{\x^0}(r)=-\infty$, then one of the following three
statements holds:
\begin{align*}&(i) &&\lim_{r\to 0}\left(\frac{u(rQ(r)\x)}{r^2}-\Pi(u,r)(Q(r)\x)\right)=
Z(\x),\\
&(ii) && \lim_{r\to 0}\left(\frac{u(rQ(r)\x)}{r^2}-\Pi(u,r)(Q(r)\x)\right)=
Z_1(\x),\\
&(iii) && \lim_{r\to 0}\left(\frac{u(rQ(r)\x)}{r^2}-\Pi(u,r)(Q(r)\x)\right)=
Z_2(\x)\end{align*}
for some $Q(r)\in \mathcal{Q}$.
However at this point we do not yet know whether
the rotation $Q(r)$ converges as $r\to 0$.

In this section we will show that in the case (ii) and (iii),
$\frac{u(r\x)}{r^2}-\Pi(u,r)(\x)$ converges
as $r\to 0$. 
In Section \ref{uniquesinglines} we will show a similar result in the case (i).

\begin{thm}\label{Z_1}
Let $n=3$ and let $u$ solve $\Delta u=-\chi_{\{u>0\}}$ in $B_1$ such that $u(0)=|\nabla u(0)|=0$
and $M:=\sup_{B_1}|u|<+\infty$.
There exist constants $r(M)>0, c(M)>0$ and $K(M)<+\infty$ such that
if
$$s\in (0,r(M)), \quad \Pi(u,s)\ge K(M) \textrm{ and } \delta^A(u(s\x))\le c(M),$$
then there is a rotation $Q$ such that either
\begin{align}
&\label{z111}\frac{u(r\x)}{r^2}-\Pi(u,r)(\x)\to Z_1(Q \x) \textrm{ as }
r\to 0
\intertext{or}
&\label{z112}\frac{u(r\x)}{r^2}-\Pi(u,r)(\x)\to Z_2(Q \x) \textrm{ as }
r\to 0.
\end{align}
Moreover, there exist $\beta>0$ and $C(M,\beta)<+\infty$ such that
$$
\sup_{B_1}\Big| \frac{\Pi(u,r)}{\sup_{B_1}|\Pi(u,r)|}-
\frac{p}{\sup_{B_1}|p|}\Big|\le 
C(M,\beta) \left(K(M)+\left|\log \left({r\over s}\right)\right|\right)^{-\beta}
$$
for all $r\in (0,s)$; here
$p(\x)=(x^2+y^2)/2-z^2$ in the case (\ref{z111}) and
$p(\x)=-(x^2+y^2)/2-z^2$ in the case (\ref{z112}).
\end{thm}
\proof 
Observe that the assumptions imply by Corollary \ref{loggrowth}
as in the proof of Corollary \ref{exist} that
$\tau_r \ge K(M)$ for $r<s$ and that
$\tau_r\to +\infty$ as $r\to 0$. Moreover we see from
Theorem \ref{symm} that
$\frac{u(s\x)}{s^2}-\Pi(u,s)(\x)$ is after rotation close to 
$Z_1(Q \x)$ or $Z_2(Q \x)$. We may assume that it is
close to $Z_1(\x)$.

We will follow the strategy
explained in the proof of Theorem \ref{symm}, and 
use the notation of that proof.
Remember that by (\ref{firststep}) and Corollary \ref{loggrowth}, up to a rotation depending on $r$,
\begin{align}\label{firststep2}
&|\Pi(u,r/2)-\tau_r p_{\delta^A_r}-c_{\delta^A_r}\big( 2p_{\delta^A_r}+\kappa(\delta^A_r) (y^2-z^2)\big)|=O(\tau_r^{-\alpha}),
\intertext{where $c_\delta\ge \bar c>0$ and} 
&\label{newkappa}
\kappa(\delta) = (1+2\delta)
\frac{3A_y(\delta)-A(\delta)}{3A_x(\delta)-A(\delta)}-
(1-2\delta)\ge \omega(\delta)>0.
\end{align}
In Theorem \ref{symm} we worked to exclude the case that $\delta^A(u(r\cdot))
\in [\beta,1/2-\beta]$ for positive $\beta$ and small $r$,
and in that $\delta$-regime, $\omega$ has been bounded from below
by a positive constant. In the present proof, however, we are interested
in the regime $\delta\to 0$, where $\omega$ degenerates.
In order to deal with this difficulty,   
we will
make a Taylor expansion of $3A_x(\delta)-A(\delta)$ and
$3A_y(\delta)-A(\delta)$ at the point $\delta=0$:
From Theorem \ref{athm} we infer that
\begin{align}\label{Taylorx}
&3A_x(\delta)-A(\delta)=3A_x(0)-A(0)-\frac{4\pi}{3\sqrt{3}}\delta+O(\delta^2)
\intertext{and}
&\label{Taylory}
3A_y(\delta)-A(\delta)=3A_y(0)-A(0)+\frac{4\pi}{3\sqrt{3}}\delta+O(\delta^2).
\end{align}
Plugging this information into
(\ref{newkappa}), we obtain that
\begin{equation}\label{kappa}
4\delta\le \kappa(\delta)\le \left(5+\frac{16\pi}{2\sqrt{3}}\right)\delta.
\end{equation}
Dividing (\ref{firststep2}) by $\tau_r$, rotating slightly and
recalling that
$p_\delta  = (1/2+\delta)x^2+(1/2-\delta)y^2-z^2$
and using that
$c_\delta\ge \bar c>0$ we infer that
\begin{equation}\label{deltadecrease}
\delta^A_{r/2}\le \delta^A_{r} - c_1 
\frac{\delta^A_{r}}{\tau_{r}}+C_2 \tau_{r}^{-1-\alpha},
\end{equation}
where $c_1>0$ is a universal constant and
$C_2<+\infty$ depends only on $M$ and $\alpha$.
Note that estimate (\ref{deltadecrease}) is independent of rotations.
In the following three Claims we will describe how (\ref{deltadecrease})
leads to a decay estimate for $\delta^A$.

\noindent\textbf{Claim 1: }{\sl There is a universal constant $\beta>0$ 
and $C_3=C_3(M)<+\infty$
such that
if $\tau_r\ge C(M)$ and $\delta^A_r\le \tau_r^{-\beta}$, then
$$
\delta^A_{r/2}\le \tau_{r/2}^{-\beta}.
$$
}

\noindent\textsl{Proof of Claim 1:} 
First, (\ref{deltadecrease}) as well as the assumption 
in the Claim imply that
$$
\delta^A_{r/2}\le \big( 1-(c_1-C_2\tau_r^{-\alpha-\beta})\tau_r^{-1} \big)\tau_r^{-\beta}.
$$
On the other hand,
$0\le \tau_{r/2}-\tau_r\le \kappa$ (see Corollary \ref{loggrowth}),
so that
\begin{equation}\label{quotientestimatefortau}
\Big( \frac{\tau_{r/2}}{\tau_r} \Big)^{-\beta}\ge 1-\beta \kappa \tau_r^{-1}.
\end{equation}
It follows that, provided that $\beta$ has been chosen small enough
(depending only on the universal constants $\kappa$ and $c_1$)
and $\tau_r$ is large enough (depending on $\kappa, c_1$ and $M$),
then
$$
\delta^A_{r/2}\le \tau_{r/2}^{-\beta},
$$
proving the claim.

$ $

Next we consider the case when $\delta^A_r\ge \tau_r^{-\beta}$.

\noindent\textbf{Claim 2:} {\sl There is a constant $\beta=\beta(\alpha)\in (0,\alpha)$ 
such that if 
$$\delta^A_r\le c_4(M),\quad \tau_r\ge C_5\textrm{ and }\delta^A_r\ge \tau_r^{-\beta},$$
then
\begin{equation}\label{buddha}
\delta^A_{r/2}\le \Big( 1- \frac{c_1}{2\tau_r}\Big)\delta^A_{r}.
\end{equation}
Moreover, if $\delta^A_{2^{-k}r}\ge \tau_{2^{-k}r}^{-\beta}$
for each $k\le k_0$ then
$$
\delta^A_{2^{-k_0}r}\le \frac{\tau_r^{2\beta}\delta^A_r}{\tau_{2^{-k_0}r}^{\beta}}
\frac{1}{\tau_{2^{-k_0}r}^\beta}.
$$}

\noindent\textsl{Proof of Claim 2:} Equation (\ref{buddha}) is a direct consequence
of equation (\ref{deltadecrease}). The last part of the Claim follows from
an induction of the first part (noting that
the assumption $\delta^A_r\le c_4(M)$ is satisfied inductively): If $\delta^A_{2^{-k}r}\ge \tau_{2^{-k}r}^{-\beta}$
for each $k\le k_0$ then 
$$
\delta^A_{2^{-k_0}r}\le \delta^A_r\Pi_{k=0}^{k_0-1}\Big(1-\frac{c_1}{2\tau_{2^{-k}r}}\Big).
$$
That product can be estimated for $\tau_r \ge C_5$, calculating
\begin{align*}
\log\Big( \Pi_{k=0}^{k_0-1}\Big(1-\frac{c_1}{2\tau_{2^{-k}r}}\Big) \Big)&=
\sum_{k=0}^{k_0-1}\log\Big(1-\frac{c_1}{2\tau_{2^{-k}r}}\Big)\le 
-\sum_{k=0}^{k_0-1}\frac{c_1}{4\tau_{2^{-k}r}}\\
&\le -\sum_{k=0}^{k_0-1}\frac{c_1}{4(\tau_{r}+\kappa k)}\le -\frac{c_1}{4}\log\Big( \frac{k_0\kappa+\tau_r}{\tau_r}\Big),\end{align*}
where we have used Corollary \ref{loggrowth}. Thus
$$
\Pi_{k=0}^{k_0-1}\Big(1-\frac{c_1}{2\tau_{2^{-k}r}}\Big)\le 
\Big( \frac{\tau_r}{k_0\kappa+\tau_r}\Big)^{c_1/4}.
$$
Choosing $\beta$ even smaller such that $2\beta\le c_1/4$
and using once more Corollary \ref{loggrowth}, we obtain
\begin{align*}
\delta^A_{2^{-k_0}r}&\le \Big( \frac{\tau_r}{k_0\kappa+\tau_r}\Big)^{2\beta}\delta^A_r\\
&\le
\Big( \frac{1}{k_0\kappa+\tau_r}\Big)^{\beta}
\frac{\tau_r^{2\beta}\delta^A_r}{(k_0\kappa+\tau_r)^{\beta}}
\le \frac{\tau_r^{2\beta}\delta^A_r}{\tau_{2^{-k_0}r}^{\beta}}
\frac{1}{\tau_{2^{-k_0}r}^\beta},
\end{align*}
proving the Claim.

$ $

\noindent\textbf{Claim 3:} {\sl There is a constant $\beta=\beta(\alpha)\in (0,\alpha)$ 
such that if $\tau_{2^{-k_0}r}\ge C(M)$ and $\delta^A_r\le c_4(M)$, then for $k\ge k_0$,
\begin{align*}
&\delta^A_{2^{-k}r}\le \tau_{2^{-k}r}^{-\beta}\textrm{ for } k \ge k_1\textrm{ and}\\
&\delta^A_{2^{-k}r}\le\frac{\tau_{r}^{2\beta}\delta^A_{r}}{\tau_{2^{-k}r}^{2\beta}}
\textrm{ for } k < k_1\\
&\textrm{for some } k_1 \le\frac{4}{\eta_0}\tau_{r}^{2}(\delta^A_{r})^{1\over \beta}.
\end{align*}}

\noindent\textsl{Proof of Claim 3:} 
We apply Claim 2 up to the first $k_1$ such that
$\delta^A_{2^{-k_1} r}\le \tau_{2^{-k_1}r}^{-\beta}$,
and we apply Claim 1 for $k\ge k_1$. 
From Claim 2 and Corollary \ref{loggrowth} we infer that
$$k_1 \le \frac{4}{\eta_0}\tau_{r}^{2}(\delta^A_{r})^{1\over \beta}.$$
Observing that the assumptions for $\tau_{2^{-k}r}$
are satisfied for $k\ge k_0$ by Corollary \ref{loggrowth}
finishes the proof of Claim 3.

$ $

In the last part of our proof we will use the decay estimate in Claim 3
in order to estimate how much $\Pi(u(r\cdot))$ moves
when varying $r$.
Let $r_k := 2^{-k}s, \tau_k := \tau_{2^{-k}s}$ and
$\delta^A_k := \delta^A_{2^{-k}s}$.
First, we infer from (\ref{firststep2}) that up to a rotation depending on $k$,
\begin{align}\label{rotationestimate}
&\sup_{B_1}\left| \frac{\Pi(u,r_k)}{\tau_k}-
\frac{\Pi(u,r_{k+1})}{\tau_{k+1}} \right|\\
&\le \sup_{B_1}\left| \frac{\tau_kp_{\delta^A_k}}{\tau_k}-
\frac{\tau_kp_{\delta^A_k}+c_{\delta^A_k}\big( 2p_{\delta^A_k}+\kappa(\delta^A_k) (y^2-z^2)\big)}{\tau_{k+1}}
\right|+
C_1(M,\alpha)\tau_k^{-(1+\alpha)},\non
\end{align}
where  $0\le c_{\delta^A_k}=-{\bf K_0}(3A_x(\delta^A_k)-A(\delta^A_k))/ (1+2\delta^A_k)
\le C_6$,
$4\delta^A_k\le \kappa(\delta^A_k)\le C_7\delta^A_k$
and $C_6,C_7$ are universal constants.
Another fact we infer from (\ref{firststep2}) is that
\begin{equation}\label{tauexpand}
\tau_{k+1} = \tau_k + 2c_{\delta^A_k} + O(\delta^A_k).
\end{equation} 
Plugging (\ref{tauexpand}) into
(\ref{rotationestimate}) yields
\begin{align*}
&\sup_{B_1}\left| \frac{\Pi(u,r_k)}{\tau_k}-
\frac{\Pi(u,r_{k+1})}{\tau_{k+1}} \right|\\
&\le
\sup_{B_1}\left| \frac{-4c^2_{\delta^A_k}+\tau_k O(\delta^A_k)}{\tau_k(\tau_k + 2c_{\delta^A_k} + O(\delta^A_k))}\right|+
C_1(M,\alpha)\tau_k^{-(1+\alpha)}
\le 
C_8(M,\alpha)\tau_k^{-(1+\alpha)}
+ C_9\frac{\delta^A_k}{\tau_k},
\end{align*}
where $C_9$ is a universal constant.
Iterating this estimate we obtain
\begin{equation}\label{someintermediateestimate}
\sup_{B_1}\left| \frac{\Pi(u,r_k)}{\tau_k}-
\frac{\Pi(u,r_{k+m})}{\tau_{k+m}} \right|
\le \sum_{i=k}^{k+m}
\frac{C_8(M,\alpha)}{\tau_{i}^{1+\alpha}}+\sum_{i=k}^{k+m}C_9
\frac{\delta^A_i}{\tau_{i}}.
\end{equation}
From Claim 3 (applied twice) and Corollary \ref{loggrowth} we conclude that,
choosing $r(M)$ small enough such that $\tau_{r(M)}\ge C(M)$,
setting $p=x^2+y^2-2z^2$ and letting $m_j\to \infty$,
\begin{align*}&\sup_{B_1}\left| \frac{\Pi(u,r_k)}{\tau_k}-
p\right|
\le C_{10}(M,\beta)\sum_{i=k}^{\infty}
\tau_{i}^{-1-\beta}+
C_9 \tau_{s}^{2\beta}\delta^A_{s}\sum_{i=k}^{k}\tau_i^{-1-2\beta}\\
&\le C_{11}(M,\beta)\sum_{i=k}^{\infty}
(M+i\eta_0/2)^{-1-\beta}
+ C_{12} \tau_{s}^{2\beta}\delta^A_{s}\sum_{i=k}^{\infty}
(M+i\eta_0/2)^{-1-2\beta}\\
&\le C_{13}(M,\beta) \tau_s^{-\beta}\textrm{ for all }k\ge k(M).
\end{align*}
Using once more Corollary \ref{loggrowth} we obtain the estimate
of the Theorem as well as
$$
\lim_{r\to 0}\left(\frac{u(\x^0+r\x)}{r^2}-\Pi(u,r,\x^0)(\x)\right)=Z_1(\x)
\textrm{ for } \x^0=0.
$$
\qed
\begin{cor}\label{twolipmanif}
Let $n=3$, let $u$ solve (\ref{main}) in $B_1$ and suppose that
\begin{equation}\label{cond}\lim_{r\to 0}\frac{u(r\x)}{\sup_{B_r}|u|} = 
(x^2+y^2)/2-z^2 \textrm{ or }
\lim_{r\to 0}\frac{u(r\x)}{\sup_{B_r}|u|} = 
z^2-(x^2+y^2)/2.\end{equation}
Then there exists an $r_0=r_0(u)$ and $f,g\in C^{0,1}(B_{r_0}')$ such that
$$
B_{r_0}\cap \{ u=0\}=\big( \{(x,y,f(x,y)): (x,y)\in B_{r_0}' \}\cup
\{(x,y,g(x,y)): (x,y)\in B_{r_0}' \}\big) \cap B_{r_0}.
$$
Moreover $f(x,y)-\sqrt{x^2+y^2}/\sqrt{2}\in C^1(B_{r_0}')$ and
$g(x,y)+\sqrt{x^2+y^2}/\sqrt{2}\in C^1(B_{r_0}')$.
The Lipschitz- and $C^1$-norms corresponding to the above statements
are uniformly bounded for
solutions $v$ sufficiently close to the fixed solution $u$
in $L^\infty(B_1)$, provided that
each $v$ satisfies
$$
\lim_{r\to 0}\frac{v(\boldsymbol{\xi}^v+rQ^v\x)}{\sup_{B_r(\boldsymbol{\xi}^v)}|v|} = 
\lim_{r\to 0}\frac{u(r\x)}{\sup_{B_r}|u|}$$
for some rotation $Q^v$ at a singular point $\boldsymbol{\xi}^v$ sufficiently
close to the origin.
\end{cor}
\proof We will show that
$\{ u=0\} \cap B_{r_0}^+=\{(x,y,f(x,y)):\; (x,y)\in B_{r_0}'\}\cap B_{r_0}$
for some $f\in C^{0,1}(B_{r_0}')$ and $f-\sqrt{x^2+y^2}/\sqrt{2}\in C^{1}$.
By symmetry a similar statement holds in $B_{r_0}^-$. We will also assume,
for the sake of definiteness, that
$$
v_r(\x)=\frac{u(r\x)}{\sup_{B_r}|u|} \to 
(x^2+y^2)/2-z^2 \textrm{ in }C^{1,\alpha}(\overline{B_1}).$$
By
the $C^{1,\alpha}$-convergence,
\begin{align*}
&\sup_{B_1} \Big| \frac{\partial v_r}{\partial z}+2z\Big|\le \omega(r)
\intertext{and}
&\sup_{B_1} \Big| v_r-\big( \frac{x^2+y^2}{2}-z^2\big) \Big|\le \omega(r)
\end{align*}
for some modulus of continuity $\omega(r)\to 0$ as $r\to 0$.
It follows that
$\{v_r=0\}\cap (\bar B_1\setminus B_{1/2})\subset\{(x,y,z)\in \bar B_1\setminus B_{1/2}:\; \dist(\cdot,\{x^2+y^2=2z^2\})\le \sigma(r)\}$
for some modulus of continuity $\sigma$. Therefore
$$
\frac{\partial v_r}{\partial z}\le -\frac{1}{4} \textrm{ on }
\bar B_1^+\setminus B_{1/2}.
$$
From the implicit function theorem and $C^{1,\alpha}$-regularity we infer
that $\{v_r=0\}\cap (\overline{B_1^+}\setminus B_{1/2})$ is a $C^{1,\alpha}$-graph
with bounded $C^{1,\alpha}$-norm (independent of $r$). It follows 
that
$\{u=0\}$ is the graph of a Lipschitz function $f$ in $B_{r_0}^+$ and we
only need to show that $f(x,y)-\sqrt{x^2+y^2}/\sqrt{2}\in C^1(B_{r_0})$.

We know that $f\in C^{1,\alpha}(\overline{B_{r_0}'}\setminus B_s')$ for every $s>0$
and that $f(rx,ry)/r$ is bounded in $C^{1,\alpha}(\overline{B_1'}\setminus B_{1/2}')$. Thus it is sufficient to show that
$$\lim_{(x,y)\to (0,0)}|\nabla (\sqrt{2}f(x,y)-\sqrt{x^2+y^2})|=0.$$
Let us consider any sequence $(x_j,y_j)\to 0$ and denote
$\sqrt{x_j^2+y_j^2}=r_j$. Then $f(r_jx,r_jy)/r_j$ will converge to
$\sqrt{x^2+y^2}/\sqrt{2}$ in $C^{1,\alpha}(\overline{B_1'}\setminus B_{1/2})$, implying that
$$
\Big(\nabla f(x,y)-\nabla \frac{\sqrt{x^2+y^2}}{\sqrt{2}}\Big)\Big|_{(x,y)=
(x_j,y_j)}\to 0.
$$
As the sequence $(x_j,y_j)$ is arbitrary, it follows that
$f-\sqrt{x^2+y^2}/\sqrt{2}\in C^{1}$.
The uniformity follows from the uniformity in Theorem \ref{Z_1}.\qed

\section{Stable Cones}\label{cone}
\begin{thm}
In $\mathbb{R}^3$ there exists a solution of (\ref{main}) in $B_1$
such that
\begin{equation}\label{stable}
0\le \frac{1}{2}\int_{B_1}|\nabla w|^2\le \int_{B_1}|\nabla w|^2-\int_{B_1\cap \{u=0\}}\frac{w^2}{|\nabla u|}\dh2
\end{equation}
for each $w\in W^{1,2}_0(B_1)$. Moreover $u\notin C^{1,1}(B_{1/2})$.
\end{thm}
Notice that the right-hand side in (\ref{stable}) is the second variation of the energy
$\int_{B_1} (|\nabla u|^2/2 - \max(u,0))$
of equation (\ref{main}). 
\proof 
By Corollary \ref{exist} there exists a solution $v$ of
$\Delta v=-\chi_{\{v>0\}}$ in $B_1$
such that the blow-up limit at the origin is $Z_1$.

Let
$$
u(\x):=\frac{v(s\x )}{s^2}
$$
for some small but fixed $s$.

For some large $M$ to be determined later
and sufficiently small $s$, Theorem \ref{Z_1}
together with Corollary \ref{loggrowth} yields that
\begin{equation}\label{gradientlarge}
|\nabla u(\x)|\ge \big( M+\log\big(\frac{1}{|\x|} \big)\big)|\x|
\end{equation}
on $\Gamma=B_1\cap \{ u=0\}$.

Choosing $s$ if necessary even smaller,
Corollary \ref{twolipmanif}
implies that $\Gamma$ consists of two Lipschitz graphs in $\overline{B_1}$.

Note that since the origin has zero capacity we may by a limiting
argument deduce that the second variation is well defined for all $w\in W^{1,2}_0(B_1)$.

If $w\in W^{1,2}$ then $w|_\Gamma\in W^{1/2,2}(\Gamma)$ by the trace theorem, which
is valid for our Lipschitz free boundary. Also from the trace theorem,
combined with Poincare's inequality,
we infer that
for each $w\in W^{1,2}_0(B_1)$
$$
\|w\|_{W^{1/2,2}(\Gamma)}\le C_1 \|w\|_{W^{1,2}(B_1)}\le C_2 \|\nabla w\|_{L^2(B_1)}.
$$ 

Using the Sobolev embedding, we obtain for $w\in W^{1,2}_0(B_1)$
\begin{equation}\label{L4estimate}
\|w\|_{L^4(\Gamma)}\le C_3\|w\|_{W^{1/2,2}(\Gamma)}\le C_4\|\nabla w\|_{L^2(B_1)}.
\end{equation}

Thus, using (\ref{gradientlarge}) and (\ref{L4estimate}),
\begin{align*}
&\int_{\Gamma}\frac{w^2}{|\nabla u|}\dh2\le \bigg(\int_{\Gamma}\frac{1}{|\nabla u|^2}\dh2 \bigg)^{1/2}
\bigg(\int_{\Gamma}w^4 \dh2\bigg)^{1/2}\\
&\le C_5\bigg(\int_{\Gamma}\frac{1}{\big( M+(\log(|\x|^{-1}))^2\big)|\x|^2} \dh2\bigg)^{1/2}\int_{B_1}|\nabla w|^2.
\end{align*}
On the other hand, $|\Gamma \cap \partial B_r|\le C_6 r$ by Corollary \ref{twolipmanif},
so that
$$
\bigg(\int_{\Gamma}\frac{1}{\big( M+(\log(|\x|^{-1}))^2\big)|\x|^2} \dh2 \bigg)^{1/2}\le
C_7\bigg(\int_{0}^1\frac{r}{\big( M+(\log|r|)^2\big)r^2} \bigg)^{1/2}\le \frac{C_8}{\sqrt{M}}.
$$
Choosing
$M\ge 4C_8^2C_5^2$ we arrive at
$$
\int_{\Gamma}\frac{w^2}{|\nabla u|}\dh2 \le \frac{1}{2}\int_{B_1}|\nabla w|^2.
$$
\qed

\section{Unique Tangent Cones at Unstable Codimension $2$ Singularities}\label{uniquesinglines}
In Theorem \ref{Z_1} we showed that if
$$\lim_{j\to\infty}\frac{u(r_j\x)}{\sup_{B_{r_j}}|u|} = 
(x^2+y^2)/2-z^2,$$
then the blow-up limit is unique and we obtain a quantitative convergence estimate.
In this section we will show the corresponding result in the case that
$$\lim_{j\to\infty}\frac{u(r_j\x)}{\sup_{B_{r_j}}|u|} = 
x^2-z^2.$$
This case corresponds to $\delta=1/2$ in the notation of
the previous sections. To make Taylor expansions
of $3A_x(\delta)-A(\delta)$ etc. around the point $\delta=1/2$ would be rather
clumsy. To get around that we will 
change the parametrization to $p_\delta=(1-\delta)x^2+\delta y^2-z^2$ and
use the $B_x,B_y,B_z,B$ defined in Section \ref{Fouriersec}.
\begin{thm}\label{uniqqesinglines}
Let $n=3$, let $u$ solve $\Delta u=-\chi_{\{u>0\}}$ in $B_1$ and suppose that 
$M:=\sup_{B_1}|u|<+\infty, \x^0\in B_{1/2},
u(\x^0)=|\nabla u(\x^0)|=0$ and that there exists a
sequence $r_j\to 0$ such that
$$
\lim_{j\to \infty}\left(\frac{u(\x^0+r_j\x)}{r_j^2}-\Pi(u,r_j,\x^0)(\x)\right)= Z(Q_{\x^0}\x),
$$
where $Q_{\x^0}$ is a rotation depending on the point $\x^0$.
Then the limit
$$
\lim_{r\to 0}\left(\frac{u(\x^0+r\x)}{r^2}-\Pi(u,r,\x^0)(\x)\right)=Z(Q_{\x^0}\x)
$$
exists (and is thus unique).

Moreover, for each $\gamma\in (0,1/4)$
there exist constants $r(M,\gamma)>0, c(M)>0, K(M)<+\infty$ 
and $C(M,\gamma)<+\infty$
such that
$$s\in (0,r(M,\gamma)), \quad \Pi(u,s)\ge K(M) \textrm{ and } \delta^B(u(s\x))\le c(M)$$
imply that 
$$
\left|\frac{\Pi(u,r,\x^0)(\x)}{\sup_{B_1}|\Pi(u,r,\x^0)|}-2(Q_{\x^0}\x)_1(Q_{\x^0}\x)_3 \right|\le C(M,\gamma) \left(K(M)+\left|\log \left({r\over s}\right)\right|\right)^{-\gamma}
$$
for all $r\in (0,s)$.
\end{thm}
\proof 
We may assume that $\x^0=0$ and that the rotation $Q_{\x^0}$ is such that
$$2(Q_{\x^0}\x)_1(Q_{\x^0}\x)_3=\x_1^2-\x_3^2.$$

By Corollary \ref{corgest}, up to a rotation depending on $r$,
\begin{equation}\label{Z_1_1}
\sup_{B_1}|\Pi(u,r/2)-\tau_r p_{\delta^B_r} - \Pi(Z_{p_{\delta^B_r}},1/2)|\le C(M,n,\gamma)
\tau_r^{-\gamma}
\end{equation}
for $\gamma<1/4$ and $r<r(M)$.
Following the strategy in the proof of Theorem \ref{Z_1}, we are going to
use (\ref{Z_1_1}) together with an analysis of
$\Pi(Z_{p_{\delta^B_r}},1/2)$ to
derive a decay estimate for $\delta^B(u(r\x))\ge 0$ (cf. Definition \ref{delta}) in $r$.
That decay in turn will make it possible to estimate how much $\Pi(u(r\cdot))$ moves
when decreasing $r$.
Note however that as the singularity examined in the present section
is by \cite{monneauweiss} {\em unstable}, we cannot expect to obtain the decay by a simple
iteration as in the proof of Theorem \ref{Z_1}. The ``pinning effect''
of the convergence assumption
$$
\lim_{j\to \infty}\left(\frac{u(\x^0+r_j\x)}{r_j^2}-\Pi(u,r_j,\x^0)(\x)\right)= Z(Q_{\x^0}\x),
$$
has to enter the proof.

\noindent\textbf{Claim:} {\sl For $r<s$,
$$\delta^B_r
\le \frac{2C(M,\gamma)\tau_r^{-\gamma}}{{\bf K_0}(3C_y(\delta^B_r)-C_0(\delta^B_r))}.
$$
}

\noindent\textsl{Proof of the Claim:} 
As the proof will be concluded by a continuity argument
in $r$, we assume that
$\Pi(u,r)\ge K(M)$ and $\delta^B(u(s\x))\le c(M)$.
From (\ref{Z_1_1}),
(\ref{piconst})
and Theorem \ref{bthm} we infer that in $B_1$, up to a rotation,
\begin{align}
\Pi(u,r/2)&=\tau_rp_{\delta^B_r}+\Pi(Z_{p_{\delta^B_r}},1/2)+O(\tau_r^{-\gamma})
\non\\
\label{rotationestimatesoonwithtildes}
&=\Big(\tau_r(1-\delta^B_r)+{\bf K_0}\big(
(3B_x(1/2-\delta)-B(1/2-\delta))x^2+(3B_y(1/2-\delta)
\\&\non \quad -B(1/2-\delta))y^2+
(3B_z(1/2-\delta)-B(1/2-\delta))z^2 \Big)+
O(\tau_r^{-\gamma})\\
&\non=\big(\tau_r(1-\delta^B_r)+{\bf K_0}(1-C_0(\delta^B_r)\delta^B_r+o(\delta^B_r)) \big)x^2\\
&\non\quad +
\big(\tau_r \delta^B_r +{\bf K_0}(3C_y(\delta^B_r)\delta^B_r-C_0(\delta^B_r)\delta^B_r) \big)y^2
\\
&\quad +\big( -\tau_r+{\bf K_0}(-3C_y(\delta^B_r)\delta^B_r+2C_0(\delta^B_r)\delta^B_r-1-o(\delta^B_r))\big)z^2+
O(\tau_r^{-\gamma}).\non
\end{align}

Rotating the coordinate system slightly to $\tilde x, \tilde y,\tilde z$,
we deduce from (\ref{rotationestimatesoonwithtildes}) that
the quotient of
the $\tilde{y}^2$
and the $\tilde{x}^2$ coefficient of $\Pi(u,r/2)$ is estimated from below by
$$
\frac{\delta^B_{r/2}}{1-\delta^B_{r/2}}\ge\frac{\tau_r\delta^B_r+3 {\bf K_0}C_y(\delta^B_r)\delta^B_r-{\bf K_0}C_0(\delta^B_r)\delta^B_r
-C_1\tau_r^{-\gamma}}{\tau_r(1-\delta^B_r)+{\bf K_0}-{\bf K_0}(C_0(\delta^B_r)\delta^B_r+o(\delta^B_r))+C_1\tau_r^{-\gamma}}.
$$

We maintain that for $\delta^B_r\le c_2(M,\gamma)$ and $\tau_r \ge C_3(M,\gamma)$,
\begin{equation}\label{contraass}
\frac{\delta^B_{r/2}}{1-\delta^B_{r/2}}\ge \frac{\delta^B_r}{1-\delta^B_r}
\textrm{ unless }
\delta^B_r
\le \frac{2C_1\tau_r^{-\gamma}}{{\bf K_0}(3C_y(\delta^B_r)-C_0(\delta^B_r))}:
\end{equation}
Subtracting the two quotients we end up with
\begin{align*}
&D:= \frac{\delta^B_{r/2}}{1-\delta^B_{r/2}}-\frac{\delta^B_r}{1-\delta^B_r}
\\
&\ge
S\Big[{\bf K_0}(3C_y(\delta^B_r)-C_0(\delta^B_r))\delta^B_r
\\&\quad -{\bf K_0}(3C_y(\delta^B_r)-2C_0(\delta^B_r))(\delta^B_r)^2-{\bf K_0}\delta^B_r+\delta^B_r o(\delta^B_r)-C\tau_r^{-\gamma}\Big],
\intertext{where}
&S=\left((\tau_r(1-\delta^B_r)+{\bf K_0}-{\bf K_0}(C_0(\delta^B_r)\delta^B_r+o(\delta^B_r))+C\tau_r^{-\gamma})(1-\delta^B_r)\right)^{-1}.
\end{align*}
For $D$ to be non-negative ---$\delta^B_r$ being by assumption small, 
$\tau_r$ being large and $C_y(\delta^B_r)\approx C_0(\delta^B_r)\approx |\log \delta^B_r|$ by Theorem \ref{bthm}---it is sufficient that
\begin{equation}\label{deltacond}
\delta^B_r
\ge \frac{2C_1\tau_r^{-\gamma}}{{\bf K_0}(3C_y(\delta^B_r)-C_0(\delta^B_r))}.
\end{equation}
Thus (\ref{contraass}) holds, and
$$\frac{2C_1\tau_r^{-\gamma}}{{\bf K_0}(3C_y(\delta^B_r)-C_0(\delta^B_r))}
\le \delta^B_r \le c_2(M,\gamma) \textrm{ implies that }
\frac{\delta^B_{r/2}}{1-\delta^B_{r/2}}\ge\frac{\delta^B_r}{1-\delta^B_r},
$$
which in turn implies that $\delta^B_{r/2}>\delta^B_r$.
Using Theorem \ref{athm} (i) together with Theorem \ref{bthm} (i) while observing that 
${\partial A \over {\partial\delta^B_r}}(\delta^B_r)=-{\partial B \over {\partial\delta^B_r}}({1\over 2}-\delta^B_r)$
and that ${\partial A_y \over {\partial\delta^B_r}}(\delta^B_r)=-{\partial B_y \over {\partial\delta^B_r}}({1\over 2}-\delta^B_r)$), we conclude that
$$
(3C_y(\delta^B_{r/2})-C_0(\delta^B_{r/2}))\delta^B_{r/2}\ge (3C_y(\delta^B_r)-C_0(\delta^B_r))\delta^B_r.
$$
As $\tau_{r/2}^{-\gamma}\le \tau_r^{-\gamma}$, it follows then that
\begin{equation}
\delta^B_{r/2}
\ge \frac{2C_1\tau_{r/2}^{-\gamma}}{{\bf K_0}(3C_y(\delta^B_r)-C_0(\delta^B_r))}.
\end{equation}
In this case 
$\delta^B_{2^{-k}r}\ge \delta^B_r$
for $k=1,2,\dots$.
Altogether we obtain that
$$\textrm{either } \liminf_{k\to \infty} \delta^B_{2^{-k}r}>0
\textrm{ or }   
\delta^B_r
\le \frac{2C_1\tau_r^{-\gamma}}{{\bf K_0}(3C_y(\delta^B_r)-C_0(\delta^B_r))}.$$
Since $\liminf_{k\to \infty} \delta^B_{2^{-k}r}>0$ would contradict our assumption
that $$\Pi(u,r_j)/\sup_{B_1} |\Pi(u,r_j)| \to x^2-z^2\textrm{ as } j\to\infty,$$
we have proved the Claim.

$ $ 

In the last part of our proof we will use the decay estimate in the Claim
in order to estimate how much $\Pi(u(r\cdot))$ moves
when varying $r$ in the interval $(0,s)$.
First note that the Claim and the fact that
$C_y\approx C_0$ when $\delta < < 1$ imply
that 
\begin{equation}\label{deltalogdelta}
C_y(\delta^B_r)\delta^B_r\le C_4 \tau_r^{-\gamma}.
\end{equation}
Next observe that by (\ref{rotationestimatesoonwithtildes})
and (\ref{deltalogdelta}),
\begin{equation}\label{taub}
\tau_{r/2} = \tau_r + {\bf K_0} + O(C_y(\delta^B_r)\delta^B_r)
= \tau_r + {\bf K_0} + O(\tau_r^{-\gamma}).
\end{equation}
Using (\ref{rotationestimatesoonwithtildes}) once more
along with (\ref{taub}) and (\ref{deltalogdelta})
we obtain
\begin{align*}
&\sup_{B_1}\left|\frac{\Pi(u,r)}{\sup_{B_1}|\Pi(u,r)|}-
\frac{\Pi(u,r/2)}{\sup_{B_1}|\Pi(u,r/2)|} \right|\\
&\le
\sup_{B_1}\bigg|\frac{\tau_r(1-\delta^B_r)x^2 + \delta^B_r y^2 - z^2}{\tau_r}\\&\quad -
\frac{\big(\tau_r(1-\delta^B_r)+{\bf K_0}-{\bf K_0}(C_0(\delta^B_r)\delta^B_r+o(\delta^B_r))+O(\tau_r^{-\gamma})\big)x^2}
{\tau_r+{\bf K_0}+O(\tau_r^{-\gamma})}\\
&\quad -\frac{
\big(\tau_r\delta^B_r+3 {\bf K_0}C_y(\delta^B_r)\delta^B_r-{\bf K_0}C_0(\delta^B_r)\delta^B_r
+O(\tau_r^{-\gamma})\big)y^2}
{\tau_r+{\bf K_0}+O(\tau_r^{-\gamma})}\\
&\quad -\frac{
\big(-\tau_r+{\bf K_0}(-3C_y(\delta^B_r)\delta^B_r+2C_0(\delta^B_r)\delta^B_r-1-o(\delta^B_r))\big)z^2}
{\tau_r+{\bf K_0}+O(\tau_r^{-\gamma})}
\bigg|
\le C_5\tau_r^{-(1+\gamma)}.
\end{align*}
As in the proof of Theorem \ref{Z_1}, an iteration leads to
$$
\Big|\frac{\Pi(u,2^{-k}s)}{\sup_{B_1}|\Pi(u,2^{-k}s)|}-(x^2-z^2) \Big|\le C_6\tau_{2^{-k}s}^{-\gamma},
$$
and we obtain the desired estimate as well as
$$
\lim_{r\to 0}\left(\frac{u(\x^0+r\x)}{r^2}-\Pi(u,r,\x^0)(\x)\right)=Z(Q_{\x^0}\x).
$$
\qed
\section{Structure of the Singular Set in $\mathbb{R}^3$}\label{regularityofthesingularset}
So far we have shown that if $\Delta u=-\chi_{\{u>0\}}$ in
$B_1\subset \mathbb{R}^3$ then the singular set
$S^u=\{x\in B_1:\; u(\x)=|\nabla u(\x)|=0 \textrm{ and } \lim_{r\to
0}\Phi^u_{\x}(r)=-\infty \}$ is divided into two parts
$S^u_1=\{\x\in S^u:\; \lim_{r\to 0}\frac{u(rQ\x+\x^0)}{\sup_{B_{r}}|u|}=
\frac{x^2+y^2}{2}-z^2\textrm{ for some
}Q\in \mathcal{Q}\}
\cup \{\x\in S^u:\; \lim_{r\to 0}\frac{u(rQ\x+\x^0)}{\sup_{B_{r}}|u|}=
-\big(\frac{x^2+y^2}{2}-z^2\big)\textrm{ for some
}Q\in \mathcal{Q}\}
$ and $S^u_2=\{\x\in S^u:\; \lim_{r\to
0}\frac{u(rQ\x+\x^0)}{\sup_{B_{r}}|u|}=xz\textrm{ for some }Q\in \mathcal{Q}\}$. In this
section we show that $S^u_1$ consists only of isolated points,
and that $S^u_2$ is locally contained in a $C^1$-curve. 
We also derive a compactness result for $S^u_2$.
\begin{lem}\label{isolated}
Let $n=3$, let $u$ solve (\ref{main}) and let $\x^0\in S^u_1$. Then there exists an
$r=r(u,\x^0)>0$ such that $\{ u=0\}\cap \{ \nabla u=0\}\cap B_{r}(\x^0)=\{\x^0\}$, that is
$\x^0$ is the only singular point in a small neighbourhood of $\x^0$.
For each class of solutions $v$ sufficiently close to $u$ in
$L^\infty(B_1)$, $\{ v=0\}\cap \{ \nabla v=0\}$ contains at most one
point in $B_r$.
\end{lem}
\proof Suppose towards a contradiction that
there exists a sequence of solutions
$u^j\to u$ in $L^\infty(B_1)$
as well as sequences $\{ u^j=0\}\cap \{ \nabla u^j=0\}\ni \x^j\to \x^0$ 
and $\{ u^j=0\}\cap \{ \nabla u^j=0\}\setminus \{ \x^j\}\ni \y^j\to \x^0$
as $j\to\infty$.
Let $r_j=|\y^j-\x^j|$. Then, 
passing if necessary to a subsequence,
$(\y^j-\x^j)/r_j\to \boldsymbol{\xi}\in \partial B_1$.
On the other hand, by $W^{2,p}$-regularity of the solution,
$u^j\to u$ in $L^\infty(B_1)\cap W^{2,p}(B_{1/2})$
so that the assumptions in Theorem \ref{Z_1}
are satisfied in $B_{\rho}(\x^j)$ for small $\rho$ and
sufficiently large $j$. 
Rotating each solution suitably around the origin,
we obtain that
$\x^j\in S^{u^j}_1$ and that
\begin{equation}\label{stupidproof}
\lim_{j\to \infty}\frac{u^j(r_j\x+\x^j)}{\sup_{B_{r_j}}|u^j|}= w(\x)=
\pm \big(\frac{x^2+y^2}{2}-z^2\big).
\end{equation}
By $C^{1,\alpha}$-convergence in equation (\ref{stupidproof}) it follows that
$w(\boldsymbol{\xi})=|\nabla w(\boldsymbol{\xi})|=0$ which is a contradiction since $|\boldsymbol{\xi}|=1$ and the origin
is the only point where
$w=|\nabla w|=0$.
\qed

We continue this section with a regularity result for $S^u_2$.
\begin{thm}\label{singularlines}
Let $n=3$. If $0\in S^u_2$ then there exists an $r(u)>0$ such that
$S^u_2\cap B_{r(u)}$ is
contained in a $C^1$-curve. For each class of solutions $v$ sufficiently close to $u$ in
$L^\infty(B_1)$,
the curves containing $S^v_2$ are relatively compact in $C^1(B_{r(u)})$. 
\end{thm}
\proof
Let us consider a sequence of solutions $u^j\to u$ in $L^\infty(B_1)$.
By uniform $W^{2,p}$-regularity of the solution, for sufficiently small
$s>0$
$$\Pi(u,s)\ge 2K(M) \textrm{ and } \delta^B(u(s\x))\le c(M)/2,$$
and for all sufficiently small $|\x^0|$ and all sufficiently large $j$,
$$\Pi(u^j,s,\x^0)\ge K(M) \textrm{ and } \delta^B(u(\x^0+s\x))\le c(M).$$
From Theorem \ref{uniqqesinglines} we obtain therefore that
$$ \sup_{B_1} \left| \frac{u^j(\y+r\x)}{\sup_{B_1} |u^j(\y+r\cdot)|}-p^j_{\y}\right|
\le \epsilon_1$$
for all sufficiently large $j$, all
$\y\in S^{u^j}_2\cap B_\rho,$ a rotation $Q^j_{\y}$,
$p(x_1,x_2,x_3)=2x_1 x_3$, $p^j_{\y}(\x)=p(Q^j_{\y}\x)$ and all $r\in (0,r_1)$.\\
\textbf{Uniform cone flatness:} \textsl{For each
$\epsilon>0$ there exists an $s_\epsilon>0$ such that
for sufficiently large $j$ and all $\y\in S^{u^j}_2\cap B_{\rho_1},$  
$\{u^j=0\}\cap\{|\nabla u^j|=0\}\cap B_s(\y)
\subset \{\y+sQ_{\y}(x_1,x_2,x_3):\; x_1^2+x_3^2\le \epsilon x_2^2\}$
for $s\in (0,s_\epsilon)$.}\\
\textsl{Proof of uniform cone flatness:} Suppose towards a contradiction that there
exists an $\epsilon_0>0$, a subsequence of solutions 
$$v^j=\frac{u^j(\y^j+s_jQ_{\y^j}^{-1}\cdot )}{\sup_{B_1} |u^j(\y^j+s_j\cdot)|}\to 2x_1 x_3 \textrm{ in } C^{1,\alpha}(\overline{B_1})$$
and a sequence of points $\boldsymbol{\xi}^j\to \boldsymbol{\xi}^0\in \partial B_1$ 
such that $v^j(\boldsymbol{\xi}^j)=|\nabla v^j(\boldsymbol{\xi}^j)|=0$ and
$(\xi_1^j)^2+(\xi_3^j)^2\ge \epsilon_0$.
Then
$(\xi_1^0)^2+(\xi_3^0)^2\ge \epsilon_0$,
contradicting
$0=|\nabla v^0(\boldsymbol{\xi}^0)|=2\sqrt{(\xi_1^0)^2+(\xi_3^0)^2}$
and thereby proving uniform cone flatness.

A standard consequence of the uniform cone flatness 
is that the class of curves containing $S^{u^j}_2\cap B_{\rho_2}$ is 
for large $j$ relatively compact
in $C^1$. An argument by contradiction yields the Theorem.
\qed
\\
The following corollary can be regarded as an extension of \cite[Corollary 7.2]{monneauweiss} outside a small cone (even outside a cusp) in the $y$-direction.
\begin{cor}\label{1234}
Let $n=3$ and suppose that for some solution $u$ to equation (\ref{main}),
$$\lim_{j\to\infty}\frac{u(r_j\x)}{\sup_{B_{r_j}}|u|} = 
2xz.$$
Then for each $\theta>0$ there exists an $r(u)>0$ such that
$\{ u=0\}\cap B_{r(u)}\cap \{|y|^2<\theta (|x|^2+|z|^2)\}$ consists of
two $2$-dimensional $C^1$-manifolds restricted to $B_{r(u)}\cap \{|y|^2<\theta (|x|^2+|z|^2)\}$, intersecting at right angles at the origin
in the $xz$-plane. For each class of solutions $v$ sufficiently close to $u$ in
$L^\infty(B_1)$ and
having each an $S^v_2$-point sufficiently
close to $0$,
the manifolds are relatively compact in $C^1$. 
\end{cor}
\proof The proof is similar to the proof of Corollary
\ref{twolipmanif} and left to the reader.\qed
\section{Appendix}\label{appendixa}
\noindent{\em Proof of Theorem \ref{athm}:}\\
Let us begin by proving (vi), that is
\begin{equation}\label{concavityeq}
3(A_y''(\delta)-A_x''(\delta))+2\delta (3A_y''(\delta)+3A_x''(\delta)-2A''(\delta))+2(3A_y'(\delta)+3A_x'(\delta)-2A'(\delta))>0
\end{equation}
for $\delta\in (0,1/2)$

For $\delta>0$ we have that $|\nabla p_\delta|\ne 0$ on
$\{p_\delta=0\}\cap \partial B_1$ and we may
thus differentiate $A_x$, $A_y$ and $A$. Thus
\begin{align*}
A_x'&(\delta)\\
&=8
\int_0^{\pi/2}\sin^3(\textrm{arccot}(\sqrt{1/2+\delta\cos(2\phi)}))
\cos^2(\phi)
\frac{\partial \textrm{arccot}(\sqrt{1/2+\delta\cos(2\phi)})}{\partial\delta}
d\phi\\
&=8\int_0^{\pi/2}\frac{\partial}{\partial \delta}\Big(
\frac{1}{12}\big(
\cos(3\textrm{arccot}(\sqrt{1/2+\delta\cos(2\phi)}))\\
&\quad -9\cos(\textrm{arccot}(\sqrt{1/2+\delta\cos(2\phi)}))
\big)\Big)\cos^2(\phi)d\phi\\
&=-\frac{32}{3}\int_0^{\pi/2}\frac{\partial}{\partial \delta}
\Big(
\frac{(2+\delta\cos(2\phi))\sqrt{1+2\delta\cos(2\phi)}}{(3+2\delta\cos(2\phi))^{3/2}}
\Big)\cos^2(\phi)d\phi,\\
A_y'&(\delta)\\
&=8\int_0^{\pi/2}\frac{\partial}{\partial \delta}\Big(
\frac{1}{12}\big(
\cos(3\textrm{arccot}(\sqrt{1/2+\delta\cos(2\phi)}))\\
&
\quad -9\cos(\textrm{arccot}(\sqrt{1/2+\delta\cos(2\phi)}))
\big)\Big)\sin^2(\phi)d\phi\\
&
=-\frac{32}{3}\int_0^{\pi/2}\frac{\partial}{\partial \delta}
\Big(
\frac{(2+\delta\cos(2\phi))\sqrt{1+2\delta\cos(2\phi)}}{(3+2\delta\cos(2\phi))^{3/2}}
\Big)\sin^2(\phi)d\phi
\intertext{and}
A'&(\delta)=-8\int_0^{\pi/2}\frac{\partial}{\partial \delta}
\big( \cos(\textrm{arccot}(\sqrt{1/2+\delta\cos(2\phi)})) \big)d\phi.
\end{align*}
Differentiating once more,
\begin{align*}
&A_x''(\delta)=-
\frac{32}{3}\int_0^{\pi/2}\frac{\partial^2}{\partial \delta^2}
\Big(
\frac{(2+\delta\cos(2\phi))\sqrt{1+2\delta\cos(2\phi)}}{(3+2\delta\cos(2\phi))^{3/2}}
\Big)\cos^2(\phi)d\phi,\\
&A_y''(\delta)=-
\frac{32}{3}\int_0^{\pi/2}\frac{\partial^2}{\partial \delta^2}
\Big(
\frac{(2+\delta\cos(2\phi))\sqrt{1+2\delta\cos(2\phi)}}{(3+2\delta\cos(2\phi))^{3/2}}
\Big)\sin^2(\phi)d\phi
\intertext{and}
&A''(\delta)=-8\int_0^{\pi/2}\frac{\partial^2}{\partial \delta^2}
\big( \cos(\textrm{arccot}(\sqrt{1/2+\delta\cos(2\phi)})) \big)d\phi.
\end{align*}

Using the last three identities
we may write the left hand side in equation (\ref{concavityeq})
as the sum of the following three terms (\ref{firstterm}), (\ref{secondterm}) and (\ref{thirdterm}):
\begin{align}\label{firstterm}
&3(A_y''-A_x'')=
32\int_0^{\pi/2}\frac{\partial^2}{\partial \delta^2}\Big(
\frac{(2+\delta\cos(2\phi))\sqrt{1+2\delta\cos(2\phi)}}{(3+2\delta\cos(2\phi))^{3/2}}
\Big)\cos(2\phi)d\phi,
\\
\label{secondterm}
&2\delta(3A_x''+3A_y''-2A'')=
-32\delta\int_0^{\pi/2}\frac{\partial^2}{\partial \delta^2}\Big(
\frac{\sqrt{1+2\delta\cos(2\phi)}}{(3+2\delta\cos(2\phi))^{3/2}}
\Big)d\phi
\intertext{and}
\label{thirdterm}
&2(3A_x'+3A_y'-2A')=
-32\int_0^{\pi/2}\frac{\partial}{\partial \delta}\Big(
\frac{\sqrt{1+2\delta\cos(2\phi)}}{(3+2\delta\cos(2\phi))^{3/2}}
\Big)d\phi.
\end{align}
In order to estimate (\ref{firstterm}), (\ref{secondterm}) and (\ref{thirdterm}),
we will the change of variables
$x=\delta\cos(2\phi)$
and then use Taylor expansions:
First we notice that ---using double factorials
$(2k+1)!!=3\times 5 \times \dots \times (2k-1) \times (2k+1)$
---
$$
\sqrt{1+2x}=1+x+\sum_{k=2}^\infty(-1)^{k+1}
\frac{(2k-3)!!}{k!}x^k
$$
and
$$
\frac{1}{(3+2x)^{3/2}}=\frac{1}{3\sqrt{3}}\sum_{k=0}^\infty (-1)^k
\frac{(2k+1)!!}{3^kk!}x^k,
$$
both sums are absolutely convergent for $|x|<1/2$.
Thus
\begin{align*}
\frac{\sqrt{1+2x}}{(3+2x)^{3/2}}&=\frac{1+x}{3\sqrt{3}}\Big(
\sum_{k=0}^\infty (-1)^k
\frac{(2k+1)!!}{3^kk!}x^k
\Big)
\\
&\quad -\Big( \sum_{k=2}^\infty(-1)^{k}
\frac{(2k-3)!!}{k!}x^k \Big)\Big(
\frac{1}{3\sqrt{3}}\sum_{k=0}^\infty (-1)^k
\frac{(2k+1)!!}{3^kk!}x^k
\Big)
\\
&=\frac{1}{3\sqrt{3}}+\frac{1}{3\sqrt{3}}\sum_{k=1}^\infty (-1)^{k+1}\Big(
\frac{(2k-1)!!}{3^kk!}(k-1)
\Big)x^k
\\
&\quad -\Big( \sum_{k=2}^\infty(-1)^{k}
\frac{(2k-3)!!}{k!}x^k \Big)\Big(
\frac{1}{3\sqrt{3}}\sum_{k=0}^\infty (-1)^k
\frac{(2k+1)!!}{3^kk!}x^k
\Big).
\end{align*}
Notice that the product of the last two sums will equal a sum
$\sum_{k=0}^\infty a_kx^k$,
where $a_k\ge 0$ when $k$ is even and $a_k\le 0$ when $k$ is odd.
Inserting $x=\delta\cos(2\phi)$ and these Taylor expansions in equation
(\ref{secondterm}) and using that $a_k\ge 0$ for even $k$
and that $\int_0^{\pi/2}a_k\cos^k(2\phi)d\phi=0$ for odd $k$,
we see that
\begin{align*}
&2\delta(3A_x''+3A_y''-2A'')
\\&=\frac{32}{3\sqrt{3}}\delta
\int_0^{\pi/2}\frac{\partial^2}{\partial \delta^2}
\Big(
\sum_{k=2}^\infty (-1)^{k}\Big(
\frac{(2k-1)!!}{3^kk!}(k-1)
\Big)\delta^k\cos^{k}(2\phi)
\Big)d\phi\\
&\quad +\frac{32}{3\sqrt{3}}\delta\int_0^{\pi/2}\frac{\partial^2}{\partial \delta^2}\Big(
\sum_{k=0}^\infty a_k\delta^k\cos^k(2\phi)
\Big)d\phi>0.
\end{align*}

Similarly we may estimate the left-hand side in (\ref{thirdterm})
and obtain that
$$
2(3A_x'+3A_y'-2A')> 0.
$$
Next, we make a Taylor expansion of the integrand in equation
(\ref{firstterm}). First, we calculate
\begin{align*}
(2+x)\sqrt{1+2x}&=(2+x)\Big(
1+x+\sum_{k=2}^\infty(-1)^{k+1}
\frac{(2k-3)!!}{k!}x^k
\Big)
\\&=2+3x+
\sum_{k=3}^\infty (-1)^{k+1}\Big(
\frac{(2k-5)!!}{k!}(3k-6)
\Big)x^k.
\end{align*}
Therefore
\begin{align*}
&\frac{(2+x)\sqrt{1+2x}}{(3+2x)^{3/2}}
\\&=\frac{1}{3\sqrt{3}}(2+3x)\Big(
\sum_{k=0}^\infty (-1)^k \frac{(2k+1)!!}{3^kk!}x^k
\Big)
\\&\quad -\Big(
\sum_{k=3}^\infty(-1)^k
\frac{(2k-5)!!}{k!}(3k-6)x^k
\Big)\Big(
\frac{1}{3\sqrt{3}}\sum_{k=0}^\infty (-1)^k
\frac{(2k+1)!!}{3^kk!}x^k
\Big)
\\&=\frac{1}{3\sqrt{3}}\Big(
2+x+\sum_{k=2}^\infty (-1)^{k+1}\Big(
\frac{(2k-3)!!}{3^kk!}(10k^2-9k+2)
\Big)x^k
\Big)
\\&\quad -\Big(
\sum_{k=3}^\infty(-1)^k
\frac{(2k-5)!!}{k!}(3k-6)x^k
\Big)\Big(
\frac{1}{3\sqrt{3}}\sum_{k=0}^\infty (-1)^k
\frac{(2k+1)!!}{3^kk!}x^k
\Big).
\end{align*}
As before, we notice that the product of the last two sums may be
written as $\sum_{k=3}^\infty a_k x^k$ where $a_k\ge 0$ when $k$ is even
and $a_k\le 0$ when $k$ is odd. Using this 
together with the above Taylor expansions,
the fact that $\int_0^{\pi/2}a_k\cos^{k+1}(2\phi)d\phi=0$
for even $k$ and that $$
b_k=\Big(
\frac{(2k-3)!!}{3^kk!}(4k^2+k-2)
\Big)>0,
$$
we obtain that the left-hand side of (\ref{firstterm}),
\begin{align*}
3(A_y''&-A_x'')
\\&=\frac{32}{\sqrt{3}}\int_0^{\pi/2}\frac{\partial^2}{\partial \delta^2}\Big(
\Big(
\sum_{k=2}^\infty (-1)^{k+1}b_k\delta^k\cos^k(2\phi)
\Big)\Big)\cos(2\phi)d\phi
\\&\quad -32\int_0^{\pi/2}\frac{\partial^2}{\partial \delta^2}\Big(
\Big(
\sum_{k=2}^\infty a_k \delta^k \cos^{k+1}(2\phi)\Big)\Big)d\phi> 0.
\end{align*}
Therefore
$$
3(A_y''-A_x'')+2\delta (3A_y''+3A_x''-2A'')+2(3A_y'+3A_x'-2A')> 0
$$
and (vi) holds.

Estimate (vii) follows now in a straightforward way:
Since $3A_x(\delta)-A(\delta)<0$ for $\delta\in (0,1/2)$
it is sufficient to show that
$$
F(\delta):= (1+2\delta)\big(3A_y(\delta)-A(\delta)\big)-(1-2\delta)\big( 3A_x(\delta)-A(\delta)\big)<0.
$$

By symmetry we have $F(0)=0$. Moreover, rotating
$\Pi(Z,1/2)=(\log(2)/\pi)xz$
(Lemma
\ref{zcalculations}) in the $xz$-plane by $45^\circ$,
we obtain that $3A_y(1/2)-A(1/2)=0$, proving (v) as well as
$F(1/2)=0$. By (vi) we also know that $F$ is
convex:
$$
F(\delta)<(1-2\delta)F(0)+2\delta F(1/2)=0,
$$
and (vii) follows.

Next we prove (i), that is
$$
\frac{\partial \big( 3A_y(\delta)-A(\delta)\big)}{\partial \delta}>0 \textrm{ for }\delta\in (1,1/2).
$$
First,
\begin{align*}
&\frac{\partial \big(3A_y(\delta)-A(\delta)\big)}{\partial \delta}\\
&=8\int_0^{\pi/2}\frac{\partial}{\partial \delta}\bigg[ 
-\frac{\sqrt{1+2\delta\cos(2\phi)}}{\big(3+2\delta\cos(2\phi) \big)^{3/2}}+
2\frac{(2+\delta\cos(2\phi))\sqrt{1+2\delta\cos(2\phi)}}{\big(3+2\delta\cos(2\phi) \big)^{3/2}}\cos(2\phi)
\bigg]d\phi.
\end{align*}
As before we may write the integrand as
$$
\frac{\partial}{\partial \delta}\bigg[-\frac{\sqrt{1+2x}}{(3+2x)^{3/2}}+
2\frac{(2+x)\sqrt{1+2x}}{(3+2x)^{3/2}}\cos(2\phi) \bigg],
$$
where $x=\delta\cos(2\phi)$. Using the Taylor series expansions calculated 
before it is easy to see that
\begin{align*}
&8\int_0^{\pi/2}\frac{\partial}{\partial \delta}\bigg[ 
-\frac{\sqrt{1+2\delta\cos(2\phi)}}{\big(3+2\delta\cos(2\phi) \big)^{3/2}}+
2\frac{(2+\delta\cos(2\phi))\sqrt{1+2\delta\cos(2\phi)}}{\big(3+2\delta\cos(2\phi) \big)^{3/2}}\cos(2\phi)
\bigg]d\phi\\
&=\int_0^{\pi/2}\bigg[ \sum_{k=0}^\infty (-1)^k a_k \delta^{k-1}\cos^k(2\phi)+
\sum_{k=0}^\infty (-1)^{k+1}b_k\delta^{k-1}\cos^{k+1}(2\phi)\bigg]d\phi,
\end{align*}
where $a_k,b_k\ge 0$. Using $\int_0^{\pi/2}\cos^j(2\phi)d\phi= 0$ for odd $j$
implies that
$$
\frac{\partial \big(3A_y(\delta)-A(\delta)\big)}{\partial \delta}>0.
$$

We argue similarly to show (ii), that is
$$
\frac{\partial \big( 3A_x(\delta)-A(\delta)\big)}{\partial \delta}<0.
$$
Here
\begin{align}
&\non\frac{\partial \big( 3A_x(\delta)-A(\delta)\big)}{\partial \delta}
\\&=-8\int_0^{\pi/2}\frac{\partial }{\partial\delta}\bigg[\frac{(1+4\cos(2\phi)+2
\cos^2(2\phi)\delta)\sqrt{1+2\delta\cos(2\phi)}}{(3+2\delta\cos(2\phi))^{3/2}}
\bigg]d\phi\non
\\&\label{3AxminusA}
=-8\int_{0}^{\pi/2}\frac{2\cos^2(2\phi)(3-2\delta)}{\sqrt{1+2\delta\cos(2\phi)}(3+2\delta\cos(2\phi))^{5/2}}d\phi.
\end{align}
Substituting the Taylor expansions
\begin{align*}
&\frac{1}{\sqrt{1+2\delta\cos(2\phi)}}=\sum_{k=0}^\infty(-1)^k\frac{(2k-1)!!}{k!}x^k
\intertext{and}
&\frac{1}{(3+2\cos(2\phi))^{5/2}}=\frac{1}{3\sqrt{3}}\sum_{k=0}^\infty
(-1)^k\frac{(2k+3)!!}{3^kk!}x^k
\end{align*}
into (\ref{3AxminusA}) we may deduce that 
\begin{align*}
&\frac{\partial \big( 3A_x(\delta)-A(\delta)\big)}{\partial \delta}
\\&=-8\int_0^{\pi/2}\frac{2\cos^2(2\phi)(3-2\delta)}{3\sqrt{3}}
\bigg[\Big(\sum_{k=0}^\infty(-1)^k\frac{(2k-1)!!}{k!}\delta^k\cos^k(2\phi) \Big)
\\& \quad \times
\Big(\sum_{k=0}^\infty
(-1)^k\frac{(2k+3)!!}{3^kk!}\delta^k\cos^k(2\phi) \Big)  \bigg]d\phi.
\end{align*}
Using that $\delta< 1/2$ and that all the odd terms in the product of the sums
equal zero we see that this expression is negative, proving (ii).

From the fact that $A_x(0)=A_y(0)$ (the projection $\Pi$ preserves symmetry) 
we infer now that
$$3A_x(0)=3/2 A_x(0) + 3/2 A_y(0)
= A_x(0) + A_y(0) - \int_{\partial B_1\cap \{ x^2+y^2> 2z^2\}} (x^2+y^2)/2$$
$$< A_x(0) + A_y(0) + \int_{\partial B_1\cap \{ x^2+y^2> 2z^2\}} -z^2
= \int_{\partial B_1\cap \{ x^2+y^2> 2z^2\}} -1
= A(0),$$
which proves (iii).

Combining (ii) and (iii) we obtain
$3A_x(\delta)-A(\delta)<0$ for $\delta\in (0,1/2)$, namely (iv).

Last, we verify (viii)-(ix):

$$
\frac{\partial }{\partial \delta}(3A_x(\delta)-A(\delta))\Big|_{\delta=0}=-
\int_0^{\pi/2}\frac{32\cos^2(\phi)\cos(2\phi)}{3\sqrt{3}}-
\frac{16\cos(2\phi)}{3\sqrt{3}}d\phi=-\frac{4\pi}{3\sqrt{3}}
$$
and
$$
\frac{\partial }{\partial \delta}(3A_y(\delta)-A(\delta))\Big|_{\delta=0}=-
\int_0^{\pi/2}\frac{32\sin^2(\phi)\cos(2\phi)}{3\sqrt{3}}-
\frac{16\cos(2\phi)}{3\sqrt{3}}d\phi=\frac{4\pi}{3\sqrt{3}}.
$$
\qed\\
{\em Proof of Lemma \ref{etalem}:}\\
From Theorem \ref{athm} we deduce that
$$
\begin{array}{l}
3A_x(1/2)-A(1/2)<0, \\
3A_y(1/2)-A(1/2)=0, \textrm{ and}\\
3A_x(0)-A(0)=3A_y(0)-A(0)<0.
\end{array}
$$
Observe now that we have proved in (\ref{secondterm}) that
$$
3A_x''+3A_y''-2A''=-
16\int_0^{\pi/2}\frac{\partial^2}{\partial \delta^2}\Big(
\frac{\sqrt{1+2\delta\cos(2\phi)}}{(3+2\delta\cos(2\phi))^{3/2}}
\Big)d\phi>0
\textrm{ for }\delta\in (0,1/2).$$ 
It follows that
$$
-\big(3A_x(\delta)-A(\delta) \big)-\big(3A_y(\delta)-A(\delta) \big)\ge
\inf_{\delta\in (0,1/2)}[-\big(3A_x(\delta)-A(\delta) \big)-\big(3A_y(\delta)-A(\delta) \big)]$$
$=: d_0 > 0$
for $\delta\in [0,1/2]$.
Now by (\ref{piconst}),
$$
\Pi(Z_{p_\delta},1/2)=-{\bf K_0}\big( (3A_x(\delta)-A(\delta))x^2+
(3A_y(\delta)-A(\delta))y^2+(3A_z(\delta)-A(\delta))z^2 \big).
$$
As $\Pi(Z_{p_\delta},1/2)$ is harmonic and thus
\begin{align*}
&(3A_z(\delta)-A(\delta))=-\big( (3A_x(\delta)-A(\delta))+
(3A_y(\delta)-A(\delta))\big)\ge d_0,
\intertext{we obtain that}
&\sup_{B_1}\big|Cp_\delta +\Pi(Z_{p_\delta},1/2)  \big|
\ge \big|Cp_\delta +\Pi(Z_{p_\delta},1/2)  \big|(0,0,1)\\
&= \big|-C-{\bf K_0}(3A_z(\delta)-A(\delta))  \big|\ge C+{\bf K_0}d_0,
\end{align*}
and the Lemma follows with
\begin{equation}\label{eta0}
\eta_0 = {\bf K_0} d_0/2.
\end{equation}
\qed\\
{\em Proof of Theorem \ref{bthm}:}\\
We begin by proving (i): For sufficiently small $\delta>0$ we have
\begin{align}\label{splittingAy}
&-\frac{\partial B_y(\delta)}{\partial \delta}=8\frac{\partial}{\partial \delta}\int_{0}^{\pi/2}\int_{\textrm{arccot}(\sqrt{(1-\delta)\cos^2(\phi)+\delta
\sin^2(\phi))})}^{\pi/2}
\sin^3(\theta)\sin(\phi)^2d\theta d\phi
\\&\non
=4\int_{0}^{\pi/2}\frac{\sin^2(\phi)}{\big(1+(1-\delta)\cos^2(\phi)+\delta\sin^2(\phi) \big)^{5/2}}\frac{\sin^2(\phi)-\cos^2(\phi)}{\sqrt{(1-\delta) \cos^2(\phi)+\delta\sin^2(\phi)}}d\phi
\\&\non\ge -C_1 + 4\int_{\pi/4}^{\pi/2}\frac{\sin^2(\phi)}{\big(1+(1-\delta)\cos^2(\phi)+\delta\sin^2(\phi) \big)^{5/2}}\frac{\sin^2(\phi)-\cos^2(\phi)}{\sqrt{(1-\delta) \cos^2(\phi)+\delta\sin^2(\phi)}}d\phi
\\&\non\ge -C_1 + c_2 \int_{\pi/4}^{\pi/2}\frac{1}{\sqrt{(1-\delta)\cos^2(\phi)+\delta\sin^2(\phi)}}d\phi
\\&\non=
-C_1 + c_2 \int_{\pi/4}^{\pi/2}\frac{1}{\sqrt{\cos^2(\phi)-\delta\cos(2\phi)}}d\phi.
\end{align}
Next we notice that when $\phi\in (3\pi/8,\pi/2)$, $-\cos(2\phi)>1/\sqrt{2}$.
Consequently the right-hand side in estimate (\ref{splittingAy})
is estimated from below by
\begin{align*}
&-C_1 + c_2 \int_{3\pi/8}^{\pi/2}\frac{1}{\sqrt{\cos^2(\phi)+\delta/\sqrt{2}}}d\phi
\\&\ge -C_1+c_3 \int_0^{1-\sin(3\pi/8)}\frac{1}{\sqrt{2t-t^2}\sqrt{\delta/\sqrt{2}+2t-t^2}}dt
\\&\ge -C_1 + c_4\int_0^{1-\sin(3\pi/8)}\frac{1}{\sqrt{t}\sqrt{c\delta+t}}dt\ge 
-C_1 + c_5 \log {1\over {\sqrt{c\delta}}}.
\end{align*}
It follows that for sufficiently small $\delta>0$,
$$
-\frac{\partial B_y(\delta)}{\partial \delta} \ge -C_6+c_7 \log\big(\frac{1}{\delta}\big).
$$
In particular,
$$
B_y(0)-B_y(\delta)\ge \Big( -C_6+c_7 \log\big(\frac{1}{\delta}\big)\Big)\delta.
$$
A similar calculation shows that
$$
B_y(0)-B_y(\delta)\le \Big( C_8-C_9 \log\big(\frac{1}{\delta}\big)\Big)\delta,
$$
so that (i) holds.

Next we are going to prove that
\begin{align}\label{derivativDiff}
&\bigg|\frac{\partial \left(B_y(\delta)-B(\delta)\right)}{\partial \delta}\bigg|
\\&=\bigg|4\int_{0}^{\pi/2}\bigg(\frac{-2\cos^2(\phi)+\delta\cos(2\phi)}{\big(1+(1-\delta)\cos^2(\phi)+\delta\sin^2(\phi) \big)^{5/2}}
\frac{\sin^2(\phi)-\cos^2(\phi)}{\sqrt{(1-\delta) \cos^2(\phi)+\delta\sin^2(\phi)}}\bigg)d\phi\bigg|\non
\\&\le C_{10},\non
\end{align}
which will imply that
\begin{equation}\label{diffCyandC0}
|C_y(\delta)-C_0(\delta)|\le C_{11}
\end{equation}
and, when combined with (i), prove (ii) and (iv).

In order to prove the inequality in (\ref{derivativDiff}), we make the change of variables $\cos(\phi)=t$,
implying that
\begin{align*}
&\bigg|4\int_{0}^{\pi/2}\bigg(\frac{-2\cos^2(\phi)+\delta\cos(2\phi)}{\big(1+(1-\delta)\cos^2(\phi)+\delta\sin^2(\phi) \big)^{5/2}}
\frac{\sin^2(\phi)-\cos^2(\phi)}{\sqrt{(1-\delta) \cos^2(\phi)-\delta\sin^2(\phi)}}\bigg)d\phi\bigg|
\\&
=\bigg| 4\int_{0}^1\frac{-2t^2+2\delta t^2-\delta}{\big(1+t^2-2\delta t^2+\delta\big)^{5/2}}
\frac{1-2t^2}{\sqrt{t^2-2\delta t^2+\delta}}\frac{1}{\sqrt{1-t^2}}dt\bigg|
\\&
\le \bigg|-4C_{12} \int_0^1\frac{\big(\delta-2(1-\delta)t^2\big)(1-2t^2)}{\sqrt{1+t}\sqrt{1-t}\sqrt{\delta+t^2-2\delta t^2}}dt\bigg|.
\end{align*}
At $t=1$ the singularity is of order $\frac{1}{\sqrt{1-t}}$ which is integrable.
We may thus estimate 
\begin{align*}
&\bigg|\frac{\partial \left(B_y(\delta)-B(\delta)\right)}{\partial \delta}\bigg|
\\&\le
C_{13}\left(1+\bigg| \int_0^{1/2}\frac{\delta-2(1-\delta)t^2}{\sqrt{\delta+(1-2\delta)t^2}}dt\bigg|\right) \le
C_{14}\Bigg(1+\sqrt{\delta})
\\&+\bigg| \int_0^{1/2}\frac{-2(1-\delta)t^2}{\sqrt{\delta+(1-2\delta)t^2}}dt\bigg|\Bigg)
\le
C_{15}(1+\bigg| \int_0^{1/2}t\bigg|)\le C_{14}.
\end{align*}

Last, we are going to show (iii), i.e.
$$
B_x(\delta)=B_x(0)+o(\delta).
$$
First we notice that for $\delta>0$, $|\nabla p_{\delta}|\ne 0$ on
the set $\{ p_\delta=0 \}$ so that we may differentiate $B_x(\delta)$:
\begin{equation}\label{intrepforaxp}
\frac{\partial B_x(\delta)}{\partial \delta}=
4\int_0^{\pi/2}\frac{\cos(2\phi)\cos^2(\phi)}
{\sqrt{(1-\delta)\cos^2(\phi)+\delta\sin^2(\phi)}(1+(1-\delta)\cos^2(\phi)+\delta\sin^2(\phi))^{5/2}}d\phi
\end{equation}
$$=
4\int_0^{\pi/2}\Bigg(\frac{\cos(2\phi)}{(1+(1-\delta)\cos^2(\phi)+\delta\sin^2(\phi))^{5/2}} \Bigg)
\Bigg( \frac{\cos^2(\phi)}{\sqrt{(1-\delta)\cos^2(\phi)+\delta\sin^2(\phi)}}\Bigg)d\phi.
$$
The term inside the first parenthesis is smooth for all $\delta$ and thus harmless. 
The term inside the second parenthesis
can be estimated by
$$
\frac{\cos^2(\phi)}{\sqrt{(1-\delta)\cos^2(\phi)+\delta\sin^2(\phi)}}
\le\frac{|\cos(\phi)|}{\sqrt{1-\delta}},
$$
which is bounded for $\delta \in (0,1/2)$.
Using the primitive function
$$
\int\frac{\cos(2\phi)\cos(\phi)}{(1+\cos(2\phi))^{5/2}}d\phi=-\frac{1}{2}\frac{\sqrt{1-\cos^2(\phi)}\cos^2(\phi)}{(1+\cos^2(\phi))^{3/2}}+C,
$$
we obtain
\begin{align*}
&B_x(\delta)=B_x(0)+\delta B_x'(0)+o(\delta)
\\&
=B_x(0)+4\delta\int_0^{\pi/2}\frac{\cos(2\phi)\cos(\phi)}{(1+\cos(2\phi))^{5/2}}d\phi+o(\delta)
=B_x(0)+o(\delta).
\end{align*}
\qed
\bibliographystyle{plain}
\bibliography{asw.bib}
\end{document}